\documentclass{compositio}

\usepackage{lscape}
\usepackage{amssymb,mathrsfs}
\usepackage{amsmath}
\usepackage[normalem]{ulem}
\usepackage{amssymb,amstext,titletoc,fancyhdr,color}

\theoremstyle{plain}
\newtheorem{theorem}{Theorem}[section]
\newtheorem{lemma}[theorem]{Lemma}
\newtheorem{proposition}[theorem]{Proposition}
\newtheorem{conjecture}[theorem]{Conjecture}
\newtheorem{corollary}[theorem]{Corollary}

\theoremstyle{definition}
\newtheorem{definition}[theorem]{Definition}
\newtheorem{example}[theorem]{Example}

\theoremstyle{remark}
\newtheorem{remark}[theorem]{Remark}

\numberwithin{equation}{section}

\def\Hom{\mathop{\rm Hom}}

\newcommand{\beqnn}{\begin{equation}}
\newcommand{\eeqnn}{\end{equation}}
\newcommand{\eb}{\begin{enumerate}}
\newcommand{\ee}{\end{enumerate}}
\newcommand{\bbm}{\begin{bmatrix}}
\newcommand{\ebm}{\end{bmatrix}}
\newcommand{\bpm}{\begin{pmatrix}}
\newcommand{\epm}{\end{pmatrix}}

\newcommand{\bi}{\begin{itemize}}
\newcommand{\ei}{\end{itemize}}
\newcommand{\beq}{\begin{eqnarray*}}
\newcommand{\eeq}{\end{eqnarray*}}

\definecolor{darkgreen}{rgb}{0,0.6,0.1}

\newcommand{\beqq}{\begin{eqnarray}}
\newcommand{\eeqq}{\end{eqnarray}}
\newcommand{\beqn}{\begin{eqnarray}}
\newcommand{\eeqn}{\end{eqnarray}}

\begin{document}

\title[Series Solutions of Linear ODEs on quotient $\mathcal{D}$-modules]{Series solutions of linear ODEs by Newton-Raphson  method on quotient $D$-modules
}

\author{Yik Man Chiang}
\email{machiang@ust.hk}
\address{Department of Mathematics, The Hong Kong University of Science and Technology,
Clear Water Bay, Kowloon, Hong Kong SAR}
\author{Avery Ching}
\email{Avery.Ching@warwick.ac.uk}
\address{Department of Statistics, The University of Warwick, Coventry, CV4 7AL, UK}
\author{Chiu Yin Tsang}
\email{h0347529@connect.hku.hk}
\address{Department of Mathematics, The University of Hong Kong, Pokfulam Road, Hong Kong SAR}


\thanks{The first and third authors were partially supported by the Research Grants Council of Hong Kong (No. 16325016).}

\begin{abstract}
We develop a $D-$module approach to various kinds of solutions to several classes of important differential equations by long divisions of different differential operators. The zeros of remainder maps of such long divisions are handled by an analogue of Hensel's lemma established recently from valuation theory. In particular, this explains the common origin of some classically known special function series solutions of Heun equations and usual Frobenius series solutions. Moreover, these remainder maps also generate  eigenvalue problems that lead to non-trivial factorizations  of certain generalized hypergeometric operators.
\end{abstract}

\maketitle

\vspace*{6pt}\tableofcontents  

\section{Introduction}
Given $L$, $K\in\mathcal{D}$, a way solutions of $L$ are expressed in terms of solutions of $K$ is to consider a left $\mathcal{D}-$linear map
	\begin{equation}\label{E:S-map}
	\mathcal{D}/\mathcal{D}L\stackrel{\times S}{\longrightarrow}\mathcal{D}/\mathcal{D}K.
	\end{equation}
	where $\mathcal{D}=\mathbb{C}\langle X,\, \partial\rangle$ denotes the standard Weyl algebra generated by $\partial$ and $X$. 

\subsection{Theory of $\mathcal{I}-$adic expansions solutions}

 However, such an $S\in\mathcal{D}$ rarely exists and so one looks for a series $S$ in an $\mathcal{I}-$adic completion for an appropriately chosen maximal ideal $\mathcal{I}$.
This map is well-defined when the remainder $R=0$ in a kind of simplified long division algorithm
	\begin{equation}\label{E:long-division}
		LS=QK+R
	\end{equation}
	for some $Q$, and the \textit{remainder} $\Phi(S)=R$ is of \textit{order} less than that of $K$. Thus the task of solving the differential equation $L$ in $\mathcal{D}/\mathcal{D}K$ can  be considered the search for \textit{zeros} of $\Phi$ in an appropriate domain.

We show that when $K$ is monic \textit{first order}, then we can find a \textit{zero} $S$ of $\Phi$ via an iterative algorithm, namely by a generalised Hensel's lemma (Newton-Raphson's method) due to a recent theory by Kuhlmann \cite{Kuhl}  for maps which are not necessarily polynomials. We note, however, that for the Netwon-Raphson method to work as in the classical case, one needs criteria that guarantees the convergence of the iterations
	\[
		S_{n+1}=\big[1-{\Phi^\prime}^{-1}\circ\Phi\big]\, (S_n), 
		\quad n=1,\, 2,\, 3,\cdots,
	\]
and this is achieved by the notion of the \textit{immediacy} of $\Phi$  \cite{Kuhl}. Under this general set up, we can  define  singularities of an element in $\mathcal{D}$ as the obstruction of being \textit{monic}.
But the main features of the comparison method proposed between two operators $L$ and $K$ in order to express solutions of one operator in terms of those of the other go much further than the classical theory of local monodromies of $L$ and $K$ which is beyond the scope of the current paper.

 We focus on $\partial=\frac{d}{dX}$ in this paper although the set up is completely general that it can be applied to the forward difference operator $\partial=\Delta$ or other types of difference operators. 

Take, for example, the hypergeometric operator \cite{AAR}, $L=\mathcal{H}=X(1-X) \partial^2+[c-(a+b+1)X]\, \partial -ab$ and $K=\partial$, so that for each $S\in\mathbb{C}[[X]]$, $\Phi(S)$ is the remainder when $LS$ is divided by $K$. Then $\Phi$ is immediate so that the Newton-Raphson iterative method of a zero of $\Phi$ yields the $n-$th approximation
	\[
		S_n=\sum_{k=0}^n\frac{(a)_k(b)_k}{(c)_k k!}X^k,
	\]
	whose limit $S$ is the classical ${}_2F_1(a,\, b;\, c; \, X)$ series.  Such $S$ gives a well-defined left $\mathcal{D}-$linear map
	\[
	\mathcal{D}/\mathcal{D}L\stackrel{\times S}{\longrightarrow}\mathcal{D}/\mathcal{D}\partial
	\]
	It is important to note that this equation is susceptible to different interpretations to $\partial$ and $X$. For example, let $f$ be analytic,  the usual $\mathcal{D}-$scalar multiplications $X f(x)=xf(x)$, $\partial f(x)=f^\prime (x)$  gives rise to
	\begin{equation}\label{E:0-shifts}
		X^n\cdot 1=x^n,
	\end{equation} and thus the classical power series solution to a differential equation at an ordinary point in $\mathbb{C}$.

The choice $K=X\partial-\lambda$  is suitable to handle the case when $L$ has a regular singularity
and a zero $S=\sum_ka_k X^k$ of the corresponding remainder map represents an analogue of the Frobenius series solution $x^\lambda S\cdot 1=x^\lambda\sum_k a_kx^k$ in the classical setting that converges in certain neighbourhood of $x=0$.  Such formal series reminds us of the $p$-adic integers, which are
canonically represented as series which converge in the $p$-adic topology. One
of the nice features of $p$-adic integers is the ease of solving equations by
an iterative method called Hensel's Lemma (Newton-Raphson's method). We apply a modified Hensel's lemma (see Theorem \ref{T:prederivative} below)  due to Kuhlmann \cite{Kuhl}  in
general topology rather than the standard Hensel lemma \cite{Schikhof}.	

Note that if we alter the $\mathcal{D}-$scalar multiplication into 
	\[
		\partial f(x)=\Delta f(x)=f(x+1)-f(x),
		\qquad
		Xf(x)=xf(x-1),
	\]so that
	\[
		X^n\cdot 1=x(x-1)\cdots (x-n+1),
	\]
then the series above represents a Frobenius solution expressed as a Newton polynomial type expansion for a \ linear difference equation also with a regular singularity at $x=0$ similarly defined. We illustrate this approach by an example to be given in \S\ref{S:CR}. But full details of different choices of difference operators and linear difference equations will appear in separate works.



There are two groups of results on the Heun equations in this paper:
	\begin{enumerate}
		\item series solutions to biconfluent Heun equations and diconfluent equations that converge in appropriately chosen topological spaces,
		 \item invariant subspaces/eigenvalues problems for Heun operators.
	\end{enumerate}

We are able to show that the remainder maps when the biconfluent Heun operators and diconfluent  operators are divided by suitable first order operators are immediate. Thus the corresponding Newton-Raphson's approximations $S_n$ give the correct $\mathcal{I}-$adic expansions (see below) that would eventually lead to known special function expansion solutions  known in the literature as to be below.
This accounts for the first group of results of the Heun equations. Although the immediacy of the remainder maps when the Heun operators and confluent Heun operators are divided by suitable first order operators are  not known, we still manage to describe the invariant subspaces and hence eigenvalues problems of these operators with the simplified long division algorithm \eqref{E:long-division}. Moreover, we resolve factorization problems of Heun and confluent Heun operator in terms   of corresponding generaralized hypergeometric operators as anticipated by Maier \cite{Maier} and Takemura \cite{Take4} as natural consequences of our eigenvalues problems. This accounts for the second group of results about the Heun operators.

\subsection{Special $\mathcal{I}-$adic expansions solutions}
 We investigate the existence of a multiplication map
	\[
	\mathcal{D}/\mathcal{D}\mathscr{H}\stackrel{\times S}{\longrightarrow}\mathcal{D}/\mathcal{D}\mathcal{H}.
	\]
	where the $\mathcal{D}=\mathbb{C}\langle X,\, \partial\rangle$, the $K$ and $L$ from \eqref{E:S-map} are respectively 
	$K=\mathcal{H}$ is the hypergeometric operator and the Heun operator 
	\[
			\mathscr{H}  = X(X-1)(X-a)  \partial^2 
		  +[\gamma (X-1)(X-a)  +\delta X(X-a)
			 +\epsilon X(X-1)]\partial 
			+ \alpha\beta X
	\]
when written in the Weyl-algebraic form
where $\alpha+\beta-\gamma-\delta-\epsilon+1=0$  holds. 
Erd\'elyi proposed (see \cite{Erdelyi1}) to expand its solutions in terms of a series of hypergeometric functions $\sum_{m=0}^{\infty}{C}_m {}_2F_1 (\alpha,\ \beta;\, \gamma-m;\, x)$ which can be expressed in terms of $D-$module language as the formal sum
	\[
		\sum_{m=0}^{\infty}{C}_m A^m,
	\] where $A=x\partial-\gamma$, $\mathcal{H}$ is the hypergeometric operator defined earlier, such that $[\mathcal{H},\, A]=-A\partial$. Such sums proposed by  Erd\'elyi obviously carries more global information about the solutions, see e.g., \cite{CCT2} for an in-depth discussion. Similar situation can be said for the other confluent Heun equations, such as \textit{confluent Heun equation} and its special case, namely the Prolate Spheroidal Wave operator (PSWE) \cite{ORX}  for which astronomers  \cite{Leaver} and signal processing engineers  \cite{Xiao_Rokhlin2003} alike need to obtain information about the solutions by expanding in terms of confluent hypergeometric functions. We solve expand solutions of the \textit{biconfluent Heun operator} (BHE)\footnote{The paper \cite{Schrodinger} is the earliest appearance of BHE  that we could find in which the author approximate the eigenvalues of the BHE in terms of those of the Hermite equations.}
	\[
		\mathscr{B}=(A-A^\dagger)^2(H+\alpha)+\beta(A-A^\dagger)+\gamma,
	\]
where  $\alpha$, $\beta$ and $\gamma$ are given, in a completion of $\mathcal{D}/\mathcal{D}(AA^\dagger+1)$	of the form
	\[
		S=\sum_{k=2}^{\infty}a_k {A^\dagger}^k,
	\]which is an element in $\overline{{A^\dagger}^2\mathbb{C}[[A^\dagger]]}$ where $H=\partial^2-X^2=AA^\dagger +1$ is the \textit{Hermite operator}, $A=\partial+X$ and $A^\dagger=\partial-X$ by proving the remainder map when $\mathscr{B}$ is divided by $AA^\dagger+1$ is immediate. This result matches essentially the sort of Hermite/Cylindrical Parabolic expansions solutions found in the literatures \cite[p. 213]{Ronveaux}, \cite{Ish_2017}. 

In the case of \textit{diconfluent Heun operator} (DHE)\footnote{we have chosen an equivalent but slightly different form \cite[(1.4.40), p. 142]{Ronveaux} that differs from the canonical form by a transformation \cite[\S1.4, Part C]{Ronveaux}
 and the BHE},
	\begin{equation} \label{E:DHE_-1}
		\mathscr{D} =X^2\partial^2+(-X^2+bX+c)\partial-aX+q,
	\end{equation}
its solutions in $\mathcal{D}/\mathcal{D}\mathcal{K}$ are studied
where $\mathcal{K}$ is the confluent hypergeometric operator
	\[
		S=\sum_{k=0}^\infty a_k \partial^k,
	\]We show that the corresponding remainder  map  $\Phi$ is immediate.  Hence we have essentially recovered and generalised these results \cite[\S3.3]{Ronveaux} and \cite[IV]{EF2008} to the $\mathcal{D}-$modules setting.

\subsection{Invariant subspace/eigenvalues problems}
We note that the termination of the proposed series of hypergeometric 
functions is an important special case and it reduces to an eigenvalue problem. Of course, our more general setup allows us to have different interpretations of the $\partial$ and $X$ compared to Erd\'elyi\rq{s} original consideration. Similarly, researchers from signal processing (e.g., \cite{ Xiao_Rokhlin2003},  \cite{ORX}) and general relativity (e.g., see Leaver \cite{Leaver} and \cite{EF2008}, \cite{EF2013}) 
also observed that it is more beneficial to expand the solutions to Prolate Spheroidal Wave operator (PSWE)
(see \cite[Eq. (1)]{Leaver}), 
	\[
		\mathscr{CH}=Z(Z-Z_0)\partial^2 +(D_1+D_2Z-Z^2)\partial +(D_3+D_4Z),
	\]
(see \cite[Eq. (143)]{Leaver}), in terms of confluent hypergeometric functions\\  $\sum_{L=-\infty}^{\infty}{a}_L\,  {}_1F_1 (\alpha;\, \, -D_4-L;\, x)$ which can again be expressed in our $\mathcal{D-}$modules language 
	\[
		\sum_{m=0}^{\infty}{a}_L A^L,
\](see \cite[Eq. (138)]{Leaver}), instead of a mere power series solutions, all being rephrased by our more general Weyl-algebraic language here.

Indeed, Erd\'elyi expanded the solutions to Heun equations in a number of different hypergeometric series (see  \cite{CCT2}, \cite{Erdelyi1}, \cite{MI2021}). Similarly, Leaver expanded solutions to the PSWE in terms of many different kinds of confluent hypergeometric functions, such as Whittaker functions and Coulomb Wave functions. 

Although we have not yet proved the immediacy of the remainder maps $\Phi$ when $L$ is diveded by $K$ the remainder maps above allow us to formulate and solve the invariant subspaces problems for Heun operator and confluent operator.

Then the recurrence relations of the coefficients in the various examples are
indeed examples of the Newton's iteration. All these series converge
automatically in an easy non-Archimedean metric, which can be related to
convergence in classical topology in a relatively straightforward manner.

The remainder map $\Phi$ above has another important application as an eigenvalue problem. For instance, we prove that a certain generalised hypergeometric operator can be factorised by a degenerate Heun operator by solving an appropriate eigenvalue problem of the remainder map $\Phi$. This also solves a particular case of a conjecture by Takemura \cite{Take4}. For instance, this gives, amongst others,  an alternative and perhaps conceptually more transparent proof of the following result by Miller and Paris (see \cite{Miller_Paris_2011}, \cite{Miller_Paris_2013})
	\[
				\begin{split}
				&{}_{r+2}F_{r+1}\left(\begin{array}{c}\alpha,\ \beta,\ e_1+m_1,\ \cdots ,\ e_r+m_r\\ \gamma,\ e_1,\ \cdots ,\ e_r\end{array};\ x\right)\\
			&=(1-x)^{\gamma-\alpha-\beta-n}
			{}_{n+2}F_{n+1}\left(\begin{array}{c}\gamma-\alpha-n,\ \gamma-\beta-n,\ f_1+1,\ \cdots ,\ f_n+1\\ \gamma,\ f_1,\ \cdots ,\ f_n\end{array};\ x\right),
				\end{split}
			\]where $m_j\in\mathbb{N}$, $f_k$ are related to the given $\alpha, \beta, \gamma, e_j$ in an algebraic manner.

This paper is divided into two main parts.  In Part I, we will prepare with the notion of ultrametric (the topological
content of the $p$-adic integers) in \S\ref{S:ultra}, and some standard notion of
$D$-modules in \S\ref{DE}. In \S\ref{SS:lde}--\S\ref{SS:Reg} we apply the unified setup to various classical examples of solutions including  power series, distributional solutions and Frobenius series solutions. In \S\ref{SS:DHE} and \S\ref{SS:BHE}, non-classical series expansion  solutions of doubly confluent Heun equations and Biconfluent Heun equations in terms of confluent hypergeometric functions ${}_1F_1$ and Parabolic cylindrical functions respectively are obtained, as in the previous subsections, as consequences of our remainder map theory.

Part II deals with more technical ``long divisions" than those considered in Part I.  \S\ref{HE}, \S\ref{hseries}  is devoted to the division of the Heun operator by special hypergeometric operators. 
In \S\ref{eigen} a finite dimensional eigenvalue problem of the remainder map $\Phi$ of Heun operator by hypergeometric operator is formulated. This provides  
 a partial
solution to a conjecture by Takemura.  This also  leads to alternative proofs of generalized hypergeometric functions ``identities". Finally in \S\ref{S:CHE}, we study the  ``long division" (remainder map) of confluent Heun operator by confluent hypergeometric operator in parallel to its Heun counterparts. A new factorization of certain  generalized ``confluent" hypergeometric operator in terms of degenerated confluent Heun operator is obtained. This is a ``confluent form" of Takemura's original conjecture. We conclude with derivation of examples of linear difference equations (Bessel operator) and their solutions arise naturally from different interpretation of $D-$modules of the remainder theory proposed in this article. Details of such applications to difference operators will appear elsewhere.

\part{$\mathcal{I}-$adic Valuations and Series Solutions of  Linear ODEs}

\section{Ultrametric Spaces and Immediate Maps} \label{S:ultra}

In the search of a series solution of a differential equation, we start with a certain
finite sum (possibly a sum with only one term.) which is meant to be an approximate
solution of the given equation. Then, one adds a ``next term'' in an optimal way in certain sense so that
the correction becomes a ``better'' approximate solution to the given equation. Here the
amended sum is improved in the sense that when it is inserted in the given
differential equation, more ``initial terms'' will vanish. This reminds us of the
classical $p$-adic integers because an integer is considered small if many initial terms
of its $p$-adic expansion vanish.

To solve a polynomial equation in the ring of $p$-adic integers, the elementary Newton
iteration is useful, that is the Hensel lemma from $p-$adic analysis guarantees the convergence of this iteration under
some mild conditions. In this article, we mimic this procedure so as to
solve differential equations.

However, we will encounter maps between spaces which are not as structured as the $p$-adic
integers. For instance, these maps are not polynomial maps and it is obscure what is the
meaning of the derivatives of these maps. Therefore, we need to weaken the conditions in
Hensel's lemma but retain the topological content at the same time. The theory of
immediate maps between ultrametric spaces will serve such a purpose, and it is what we
will briefly survey in this section. The readers may find the detail in
\cite{Kuhl, Schikhof}.

\begin{definition}[{\cite[\S 1.1, p. 1731]
{Kuhl} (see also \cite{Kuhlmann_book})}]
Let $\Gamma$ be a totally ordered set with a last element $\infty$.
A {\it $\Gamma$-valued ultrametric space} is a pair $(\mathfrak{X},u)$ in which $\mathfrak{X}$ is a set,
$u:\mathfrak{X}\times \mathfrak{X}\to\Gamma$ is a map such that for every $X$, $Y$, $Z\in \mathfrak{X}$,
\begin{itemize}
\item[(U1)]
$u(X,Y)=\infty$ if and only if $X=Y$;
\item[(U2)]
$u(X,Z)\geq\min\{u(X,Y),u(Y,Z)\}$;
\item[(U3)]
$u(X,Y)=u(Y,X)$.
\end{itemize}
\end{definition}
\bigskip

\begin{example}\label{Eg:valuation}
Let $R$ be a commutative ring with unity and $\mathfrak{I}\subset \mathfrak{R}$ is a maximal ideal.
Define $u:\mathfrak{R}\times \mathfrak{R}\to\mathbb{N}\cup\{\infty\}$ by
\[
u(X,Y)=\max\{n:X-Y\in \mathfrak{I}^n\}.
\]
Then $u$ is an ultrametric. In particular, when $\mathfrak{R}=\mathbb{Z}$, and $\mathfrak{I}=p\mathbb{Z}$ is the ideal generated by a prime number $p$, the $p$-adic ultrametric is recovered.
\end{example}

\begin{definition}
Let $u:\mathfrak{X}\times \mathfrak{X}\to\Gamma$ be an ultrametric. For each $X\in \mathfrak{X}$ and $\alpha\in\Gamma$,
define the closed ball
\[
B_{\alpha}(X)=\{Z\in \mathfrak{X}:u(X,Z)\geq\alpha\}.
\]
For each $X$, $Y\in \mathfrak{X}$, we also denote
\[
B(X,Y)=B(Y,X)=B_{u(X,y)}(X)=B_{u(X,Y)}(Y).
\]
\end{definition}
\bigskip

\begin{definition}[{\cite[\S1.1, p. 1732]{Kuhl}}]\label{D:immediate}
Let $(\mathfrak{X},u)$, $(\mathfrak{Y},v)$ be ultrametric spaces and let $f:\mathfrak{X}\to \mathfrak{Y}$ be a map. Then $Y\in \mathfrak{Y}$ is called
an {\it attractor} of $f$ if for each $X\in \mathfrak{X}$ such that $f(X)\neq Y$, there exists $X'\in \mathfrak{X}$
such that 
	\begin{itemize}
	\item[(AT1)] $v(f(X'),Y)>v(f(X),Y)$;
	\item[(AT2)] $f(B(X,X'))\subset B(f(X),Y)$.
	\end{itemize}
The map $f$ is said to be {\it immediate} if every $Y\in \mathfrak{Y}$ is an attractor for $f$.
\end{definition}
\bigskip

\begin{example}[(Hensel's Lemma)]
Let $p$ be a prime. The $p$-adic ultrametric is defined on $\mathbb{Z}$ as in
the last example. Suppose that $f:\mathbb{Z}\to\mathbb{Z}$ is a polynomial map.
Then, for each $x$, $y\in\mathbb{Z}$ such that $f(x)\neq y$ and $f'(x)\neq 0$,
\[
x'=x-f(x)(f'(x))^{-1}\hspace{1cm}\mbox{(Note: $(f'(x))^{-1}$ is the inverse
of $f'(x)$ mod $p$.)}
\]
is an integer satisfying (AT1) and (AT2). Thus, $y$ is an attractor for $f$ 
(see \cite[Theorem 7.3]{Eisenbud}).
\bigskip

It is natural to raise the question asking if the condition $f'(x)\neq 0$ can be
weakened, especially when $f$ is not a polynomial map so that one cannot formulate
the derivative $f'$ easily. Theorem \ref{T:prederivative} below will give a partial
answer to this question.
\end{example}
\bigskip

\begin{definition}[{\cite[\S1.1, p. 1732]{Kuhl}}]
An ultrametric space is said to be {\it spherically complete} if every nest of balls has nonempty
intersection.
\end{definition}
\bigskip

\begin{theorem}[{\cite[Theorem 12]{Kuhl}}]
Let $f:\mathfrak{X}\to \mathfrak{Y}$ be an immediate map so that $\mathfrak{X}$ is spherically complete. Then $f$ is surjective and $\mathfrak{Y}$ is spherically complete.
\end{theorem}
\bigskip

\begin{definition}
Let $G$ be an abelian group. A {\it $\Gamma$-valued valuation} on $G$ is a map $\nu:G\to\Gamma$
such that
\begin{itemize}
\item[(V1)] $\nu(X)=\infty$ if and only if $X=0$;
\item[(V2)] $\nu(X-Y)=\min\{\nu(X),\nu(Y)\}$.
\end{itemize}
A valuation always induces an ultrametric $u:G\times G\to\Gamma$ by
\[
u(X,Y)=\nu(X-Y)\mbox{ and vice versa}.
\]
In particular, an ideal $I$ in a commutative ring $R$ induces a valuation called
the $I$-adic valuation.
\end{definition}
\bigskip

\begin{theorem}[{\cite[Proposition 11]{Kuhl}}]
Let $(G,\nu)$, $(G',\nu')$ be valued abelian groups and let $f:G\to G'$ be a group homomorphism.
Suppose that for each $Y\in G'\backslash\{0\}$, there exists $X\in G$ such that
\begin{itemize}
	\item[(IH1)] $\nu'(Y-f(X))>\nu'(Y)$;
	\item[(IH2)] for every $Z\in G$, $\nu(X)\leq\nu(Z)$ implies $\nu'(f(X))\leq\nu'(f(Z))$.
\end{itemize}
Then $f$ is immediate.
\end{theorem}

In this paper, the following theorem is often used in place of the Hensel's lemma.
\bigskip

\begin{theorem}[{\cite[Lemma 13]{Kuhl}}] \label{T:prederivative}
Let $(G,\nu)$, $(G',\nu')$ be valued abelian groups and let $f$, $\overline{f}:G\to G'$ be group
homomorphisms. If $f$ is immediate and for every $X\in G$,
\[
\nu'(\overline{f}(X)-f(X))>\nu'(f(X))\mbox{ or }\overline{f}(X)=f(X)=0,
\]
then $\overline{f}$ is immediate.
\end{theorem}
\bigskip

\begin{definition}
Let $f:G\to G'$ be a map. Then $0\neq X\in G$ is said to be {\it regular} if it satisfies (IH2).
\end{definition}

\section{Differential Operators: Solutions, Orders and Division Algorithm}\label{DE}
We gather some standard definitions of $\mathcal{D}-$modulus \cite{coutinho}, \cite{kashiwara}
that we shall use in this section.
\bigskip

\begin{definition}
Let $n\in\mathbb{N}$. The {\it Weyl algebra} $\mathcal{D}_n$ is the $\mathbb{C}$-algebra
with $2n$ generators $X_1,\cdots,  X_n$, $\partial_1$,..., $\partial_n$ subject
to the relations
\[
	[X_i,X_j]=0;\hspace{1cm}[\partial_i,\partial_j]=0;\hspace{1cm}[\partial_i,X_j]=\delta_{ij}.
\]
In this paper, we focus on $n=1$ and denote $\mathcal{D}_1$  by $\mathcal{D}$. Note that the
symbols $X$ and $\partial$ can sometimes be misleading as they are only a pair of symbols
following the syntax above instead of the semantics of ``variable'' and ``derivative''.
\end{definition}
\bigskip

\begin{definition}[(\textbf{Solution})]
Let $M$, $N$ be left $\mathcal{D}_n$-modules. A {\it solution} of $M$ in $N$ is a
left $\mathcal{D}_n$-linear maps from $M$ to $N$. In particular, the set of all solutions
of $M$ in $N$ is denoted by
\[
\Hom(M,N).
\]
\end{definition}
\bigskip

\begin{example}
Let $L$, $K\in\mathcal{D}$ be nonzero ``differential operators". Suppose that there
are differential operators $S$, $Q\in\mathcal{D}$ such that
\[
LS=QK.
\]
Then the map
\[
\mathcal{D}/\mathcal{D}L\stackrel{\times S}{\to}\mathcal{D}/\mathcal{D}K
\]
is left $\mathcal{D}$-linear. For each left $\mathcal{D}$-module, it induces a map between
solution sets
\[
\Hom(\mathcal{D}/\mathcal{D}K,N)\to
\Hom(\mathcal{D}/\mathcal{D}L,N).
\]
\end{example}
\bigskip

\begin{definition}[(\textbf{Order})] \label{D:order}
An {\it order} in $\mathcal{D}$ is an increasing filtration
$F^0\mathcal{D}\subset F^1\mathcal{D}\subset F^2\mathcal{D}\subset
...\subset\mathcal{D}$ so that 
\begin{itemize}
\item[(O1)] the set of differential operators of order at most zero $F^0\mathcal{D}$ forms
a commutative $\mathbb{C}$-algebra such that
\[
F^0\mathcal{D}=\{X\in\mathcal{D}:XY=YX\mbox{ for all }Y\in F^0\mathcal{D}\};
\]
\item[(O2)] $F^i\mathcal{D}\cdot F^j\mathcal{D}\subset F^{i+j}\mathcal{D}$;
\item[(O3)] $[F^i\mathcal{D},F^j\mathcal{D}]\subset F^{i+j-1}\mathcal{D}$
\end{itemize}
hold.  An operator in $\mathcal{D}$ is said to have \textit{order at most $k$} if it is in $F^k\mathcal{D}$.
\end{definition}
\bigskip

\begin{example}[(\textbf{Standard order of a differential operator})] \label{Ex}
Let $F^0\mathcal{D}$ be the $\mathbb{C}$-subalgebra of $\mathcal{D}$ generated by
$X$. Then an operator $L\in\mathcal{D}$ is of \textit{order at most} $k$ (i.e.
in $F^k\mathcal{D}$) if
\[
[L,f]\in F^{k-1}\mathcal{D}\hspace{1cm}\mbox{for every }f\in F^0\mathcal{D},
\]
and $F^k\mathcal{D}$ from the corresponding filtration turns out to be an order defined from Definition \ref{D:order} in $\mathcal{D}$. This order of $L$ is also  called the
{\it standard/usual order}.
\end{example}
\bigskip

\begin{example} \label{Ex1}
Let $F^0\mathcal{D}$ be the $\mathbb{C}$-subalgebra of $\mathcal{D}$ generated by
$\partial$. $F^k\mathcal{D}$ are then defined as in the last
example. This defines an order in $\mathcal{D}$ in which the order of a monomial
is the highest power of $X$ involved.
\end{example}
\bigskip

\begin{definition}[(\textbf{Monic operator})]\label{D:monic}
Let $\mathcal{D}$ be endowed with an order. A first order operator $L\in\mathcal{D}$
is said to be {\it monic} if the composition
\[
F^0\mathcal{D}\longrightarrow\mathcal{D} \longrightarrow\mathcal{D}/\mathcal{D}L
\]
is a $\mathbb{C}$-vector space isomorphism.
\end{definition}
\bigskip

\begin{example}
Let $\mathcal{D}$ be equipped with the usual order. Then $\partial$, and for that matter similarly for $k\partial\ (k\in\mathbf{C}\backslash\{0\})$,  are monic but $X\partial$
is not.
\end{example}
\bigskip

By the Definition \ref{D:monic}, we obtain the following
\bigskip

\begin{proposition}[(\textbf{Simplified Division Algorithm})]\label{division}
Let $\mathcal{D}$ be equipped with an order. Given $F$, $G\in\mathcal{D}$ in which $G$ is a
monic first-order operator, then there exists unique $Q\in \mathcal{D}$, $R\in F^0\mathcal{D}$ such that
\[
F=QG+R.
\]
\end{proposition}

By Proposition \ref{division}, it makes sense to define the {\it remainder map} as follows.
\bigskip

\begin{definition}[(\textbf{Remainder Map})] \label{D:remainder}
Let $\mathcal{D}$ be equipped with an order. Fix $L$, $K\in\mathcal{D}$ in which $K$
is a monic first-order operator. For each $X\in F^0\mathcal{D}$, 
there exists unique $Q\in\mathcal{D}$ and $r\in F^0\mathcal{D}$ such that
\begin{equation} \label{Div}
LS=QK+R.
\end{equation}
The map $\Phi: F^0\mathcal{D}\to F^0\mathcal{D}$ defined by $\Phi(S)=R$ 
 is called
the {\it remainder map} of ``$L/K=L\div K$''.
\end{definition}
\bigskip

\begin{example} If $L=p_n(X)(\partial)^n+p_{n-1}(X)(\partial)^{n-1}+...+p_0(X)$ and
$K=\partial$, then 
	\[
		\Phi(X^m)=p_n(0)(m)(m-1)\cdots (m-n+1) X^{m-n}\quad \mod X^{m-n+1},
\ \mbox{for all }m\ge n\ge 1.
	\]
\end{example}
\bigskip

\begin{remark}
Refer to the definition above, $F^0\mathcal{D}$ is a commutative ring with unity.
A choice of a maximal ideal $I\subset F^0\mathcal{D}$ induces an ultrametric in
$F^0\mathcal{D}$. It is an interesting question to ask if the remainder map $\Phi$
is immediate. In case if the answer is affirmative, $\Phi$ will in particular be
surjective after it is extended to the ``completion" of $F^0\mathcal{D}$. In other words,
there exists $S$ in the completion of $F^0\mathcal{D}$ such that $LS=QK$. Such a
factorization yields a left $\mathcal{D}$-linear map
\[
\mathcal{D}/\mathcal{D}L\stackrel{\times S}{\longrightarrow}\mathcal{D}/\mathcal{D}K.
\]
Therefore, for each left $\mathcal{D}$-module $N$, solutions of $L$ in $N$ will
be understood via the induced map
\[
\Hom(\mathcal{D}/\mathcal{D}K,N)\longrightarrow
\Hom(\mathcal{D}/\mathcal{D}L,N).
\]

Several examples of this type will be presented in the next section.
\end{remark}
\bigskip

\section{$\mathcal{I}-$adic Valuation at Ordinary points  and Immediate Maps}
In the previous section, we revised the notion of $D$-modules. In particular,
quotient modules were emphasized. It becomes the unified setup of the various
spaces which accommodate the different kinds of special functions. For instance,
$\mathcal{D}/\mathcal{D}\partial$ is the space of polynomials which becomes the
space of power series after a certain completion. As another example, the
\[ \mathcal{D}/\mathcal{D}[X(1-X)\partial^2+(c-(a+b+1)X)\partial-ab] \]
is a
space containing $_2F_1(a,b;c;\, x)$ and the various actions by differential operators.
It becomes a space containing infinite series of hypergeometric functions after
completion \cite{AAR}. In this section, we demonstrate examples of finding solutions of
linear differential equations in these $D$-modules. Inevitably, one should start
with the elementary power series solutions.
\bigskip

\subsection{Power series as $X-$adic expansions  at an ordinary point and immediate maps}\label{SS:lde}

Throughout this section, let $\mathcal{D}$ be equipped with the usual order,
so that $F^0\mathcal{D}=\mathbb{C}[X]$ is the ring of polynomials in one variable
(see Example \ref{Ex}).
\bigskip

\begin{definition}[(\textbf{Ordinary point})] \label{D:p} Suppose that a $L\in\mathcal{D}$ is a differential operator of order $n$ (Definition \ref{D:order}). Then
	\begin{enumerate}
		\item we define  $p(X)=p_L(X)$ to be a polynomial $p(X)\in\mathbb{C}[X]$ so that $L-p(X)\partial^n$ is of order at
most $n-1$ in the sense of the Definition \ref{D:order}. 
		\item and  $X=0$ is said to be an \textit{ordinary point} of $L$ if $p(0)\not=0$.
	\end{enumerate}
\end{definition}

\bigskip

It is clear that the polynomial $p(X)$ defined in the above definition is  unique with respect to a given $L$.

\medskip

Now we plan to study the remainder map $\Phi:\mathbb{C}[X]\to\mathbb{C}[X]$ of
$L/\partial=L\div \partial$ in Definition \ref{D:remainder} where $K=\partial$. Here $\mathbb{C}[X]$ is also equipped with the ultrametric induced
by the ideal generated by $x$.

For each $f(X)\in\mathbb{C}[X]$, consider the remainder when $Lf$ is divided by $\partial$. It is easy to see that $p(0)\partial^n$ is the
term in $L$ which decreases the $X-$adic valuation of such remainder the most when  $p(0)\not=0$ and less so if $p(0)=0$. We formulate this observation into
the following lemma.
\bigskip

\begin{lemma}
Let $\nu$ be the $X$-adic valuation in the maximal ideal $X\mathbb{C}[X]\subset \mathbb{C}[X]$ in Example \ref{Eg:valuation}.  Let 
$\Phi$ be the remainder map $\Phi:\mathbb{C}[X]\to\mathbb{C}[X]$ of $L\div \partial$. Then for each $f(X)\in\mathbb{C}[X]$,
\[
\nu(\Phi(f(X)))\geq\nu(f(X))-n
\]with equality if and only if $x=0$  is an ordinary point of $L\in\mathcal{D}$, i.e., $p(0)=p_L(0)\not=0$.
In particular, $\Phi$ is continuous.
\end{lemma}
Hence $\Phi$ is extended to a map defined on the completion
$\overline{\mathbb{C}[X]}=\mathbb{C}[[X]]$. We call this extension $\Phi$ also.

\bigskip

\begin{lemma}[(\textbf{Tangential approximation  at ordinary point})]\label{L:tangential}  Let $L\in\mathcal{D}$ and 
$\Phi$ be the remainder map $\Phi:\mathbb{C}[X]\to\mathbb{C}[X]$ of $L\div\partial$  and $p(X)=p_L(X)$. Suppose that $x=0$ is an ordinary point of $L$, i.e.,
$p(0)\neq 0$. Define a map $\Phi':\mathbb{C}[X]\to\mathbb{C}[X]$ by
	\begin{equation}\label{E:tangent_map}
		\Phi'(f(X))=p(0)f^{(n)}(X)\hspace{1cm}\mbox{for all }f(X)\in\mathbb{C}[X].
	\end{equation}
Let $\nu$ be the $X-${adic} valuation. Then for each $f\in\mathbb{C}[X]$,
	\[
		\nu(\Phi(f(X))-\Phi'(f(X)))>\nu(\Phi(f(X)))=\nu\big(\Phi^\prime(f(X))\big).
	\]
\end{lemma}
\bigskip

\begin{proof}
For each $f(X)\in\mathbb{C}[X]$, let $f(X)=X^mg(X)$ for some $m\in\mathbb{N}$ and
$g(0)\neq 0$ (so that $\nu(f(X))=m$.). Then it is obvious that
	\[
		\begin{array}{rll}
 		&\Phi(f(X))&\\
		=&\Phi(X^mg(X))& \\
		=&p(0)m(m-1)...(m-n+1)x^{m-n}g(0)&\mod X^{m-n+1}\\
	=&\Phi'(f(X))&\mod X^{m-n+1}.
	\end{array}
	\]
Thus, 
		\begin{align*}
			\nu\big(\Phi(f(X))-\Phi'(f(X))\big) &\geq m-n+1\\
			&>m-n\\
			&=\nu\big(\Phi(f(X))\big)=\nu\big(\Phi^\prime(f(X))\big).
		\end{align*}
\end{proof}

\bigskip

\begin{definition}[(\textbf{Tangential map at regular singularity})] We call the map $\Phi':\mathbb{C}[X]\to\mathbb{C}[X]$ defined in (\ref{E:tangent_map}) above a \textit{tangential map} of $\Phi$.
\end{definition}
\bigskip

The $\Phi'$ therefore serves as a \textit{first-order approximation} to $\Phi$.
\bigskip

\begin{lemma}[(\textbf{Tangential map is immediate})]\label{L:immediate2}
 Suppose $L\in\mathcal{D}$ has an ordinary point at $X=0$ (i.e., $p(0)\neq 0$).  Then the tangential map $\Phi':\mathbb{C}[X]\to\mathbb{C}[X]$ (\ref{E:tangent_map}) defined in the last Lemma 
is immediate.
\end{lemma}
\bigskip

\begin{proof}

For each $f(X)\in\mathbb{C}[X]$, we show that $f(X)$ is an attractor of $\Phi'$.\\
	\begin{itemize}
		\item[(AT1)] For each $G(X)\in\mathbb{C}[X]$ such that $\Phi^\prime (G(X))\neq f(X)$, let
$F(X)$ be an $n$-folded antiderivative of $f(X)/p(0)$ such that
$\nu(F(X)-G(X)):=m\geq n$. Then,
	\[
		\nu(\Phi'(F(X))-f(X))=\nu(0)=
		\infty>\nu(\Phi'(G(X))-f(X)).
	\]
		\item[\qquad (AT2)] For each $H(X)\in B(F(X),G(X))$, write
		\begin{align*}
			G(X)&=F(X)+X^mg(X)\\
			H(X)&=F(X)+X^mh(X)
		\end{align*}
		for the same $m\geq n$ and $g(X),\, h(X)\in\mathbb{C}[X],\  g(0)\neq 0$.
Then,
	\[
		\nu(\Phi'(H(X))-f(X))\geq m-n=\nu(\Phi^\prime(G(X))-f(X)).
	\]
	\end{itemize}
Thus, $\Phi^\prime(B(F(X),G(X)))\subset B(\Phi^\prime(G(X)),f(X))$. Hence $\Phi^\prime$ is immediate.
\end{proof}
\bigskip

Recall that the remainder map $\Phi$ can be extended to a map $\tilde{\Phi}$ defined on the completion
$\overline{\mathbb{C}[X]}=\mathbb{C}[[X]]$. We continue to call this extension with the same notation $\Phi$ in the following discussion.
\bigskip

\begin{theorem}[(\textbf{Remainder map is immediate})] Let $L\in\mathcal{D}$. If $X=0$ is an ordinary point (i.e., $p(0)\neq 0$) of $L$.  Then the (extended) remainder map $\Phi:\mathbb{C}[[X]]\to\mathbb{C}[[X]]$ is immediate and hence surjective.
\end{theorem}
\bigskip

\begin{proof} This is a direct consequence of the Theorem \ref{T:prederivative},  Lemma \ref{L:tangential} and Lemma \ref{L:immediate2}.
\end{proof}
\bigskip

As a direct consequence, if $p(0)\neq 0$, then for each polynomial $q(X)$ of degree at most $n-1$, hence
there exists a non-zero $S\in\mathbb{C}[[X]]$ such that $\nu(S(X)-q(X))\geq n$ and 
$\Phi(S)=0$. Such $S$ induces a left $\mathcal{D}$-linear map
	\[
		\overline{\mathcal{D}/\mathcal{D}L}\stackrel{\times S}{\longrightarrow}
		\overline{\mathcal{D}/\mathcal{D}\partial},
	\]and thus it is a ``power series solution'' of $L$ that we encountered in analysis.

Note also that we have a classical theorem concerning how $S$ is described as a limit of polynomials:
\bigskip

\begin{theorem}[(\textbf{Radius of convergence})]
Let $S=\sum_{k=0}^{\infty}a_kX^k$ be the $X-${adic} expansion of an element in the completion $\mathbb{C}[[X]]$ so that
	\[
		\overline{\mathcal{D}/\mathcal{D}L}\stackrel{\times s}		{\longrightarrow}
		\overline{\mathcal{D}/\mathcal{D}\partial}
	\]
is left $\mathcal{D}$-linear. Then
	\[
		\limsup_k|a_k|^{1/k}\leq 1/R
	\]
where $R=\min\{|a|:p(a)=0\}$.
\end{theorem}
\bigskip

\begin{proof}
Assume on the contrary, $\limsup_k|a_k|^{1/k}=1/r>1/R$ so that the series $s$ gives
an analytic solution $f$ of $L$ which is defined on $D(0,r)$. Now for each
$w\in\partial D(0,r)$, the differential operator $L$ has a regular point at $w$.
So $L$ has a complete set of analytic solutions around $w$. In other words,
there exists $\delta_w>0$ and an analytic function
$f_w:D(w,\delta_w)\to\mathbb{C}$ such that
$f|_{D(0,r)\cap D(w,\delta_w)}=f_w|_{D(0,r)\cap D(w,\delta_w)}$. But then
$\overline{D(0,r)}$ is compact and
\[\mathbb{C}\backslash (\cup_w D(w,\delta_w)\cup D(0,r))\] 
is closed.
Hence the distance between these two sets is a positive number
$\delta>0$. We see that $f$ can be extended to $D(0,r+\frac{\delta}{2})$,
contrary to the fact that the radius of convergence of the series $s$ is $r$.
\end{proof}
\bigskip

\subsection{Application to Airy operator and immediate maps}\label{SS:airy}
\bigskip
The Weyl algebraic form of Airy's differential equation
\[
    y'' - xy = 0
\]
is represented by the symbol
\[
    L = \partial^2 - X.
\]
Let $\Phi:\mathbb{C}[[X]]\to X\mathbb{C}[[X]]$ be the remainder map of $L/\partial$.\ \ Since
\[
    LX^m = \partial^2 X^m - X^{m+1} \equiv m(m-1)X^{m-2} - X^{m+1}\qquad \mbox{ mod } \partial,
\]
we have $\Phi(X^m) = m(m-1)X^{m-2} - X^{m+1}$.\ \ The tangential map $\Phi':\mathbb{C}[X]\to X\mathbb{C}[X]$ is the linear map defined on the standard basis by
\[
    \Phi'(X^m) = m(m-1)X^{m-2}.
\]
Hence an inverse of $\Phi^\prime$ can be given by
\[
    \Phi'^{-1}(X^m) = \frac{1}{(m+2)(m+1)}X^{m+2}.
\]
\begin{enumerate}
\item Now choose an ``initial trial" $S_0=1$.\ \ Then $\Phi(S_0)=-X$ and so
\[
    S_1 = (1-\Phi'^{-1}\circ\Phi)S_0 = 1 - \Phi'^{-1}(-X) = 1 + \frac{1}{3\cdot 2}X^3.
\]
Since
\[
    LS_1 = L + \frac{1}{3\cdot 2}LX^3 \equiv -X + \frac{1}{3\cdot 2}(3\cdot 2 X - X^4) \equiv - \frac{1}{3\cdot 2}X^4 \mod \partial,
\]
we have $\Phi(S_1)=-\frac{1}{3\cdot 2}X^4$ and so
\[
    S_2 = (1-\Phi'^{-1}\circ\Phi)S_1 = \left(1 + \frac{1}{3\cdot 2}X^3\right) - \Phi'^{-1}\left(-\frac{1}{3\cdot 2}X^4\right) = 1 + \frac{1}{3\cdot 2}X^3 + \frac{1}{(3\cdot 2)(6\cdot 5)}X^6.
\]
Inductively, we obtain an $n-$th Newton\rq{s} approximation to be
\[
    S_n = \sum_{k=0}^n{\frac{3^k(\frac{1}{3})_k}{(3k)!}X^{3k}}.
\]
From this we obtain the first well-known power series solution
\[
    Y_1(X) = \sum_{k=0}^\infty{\frac{3^k(\frac{1}{3})_k}{(3k)!}X^{3k}}
\]to the Airy operator.
\item We next choose an ``initial trial" $S_0=X$.\ \ Then $\Phi(S_0)=-X^2$ and so
\[
    S_1 = (1-\Phi'^{-1}\circ\Phi)s_0 = x - \Phi'^{-1}(-X^2) = x + \frac{1}{4\cdot 3}X^4.
\]
Since
\[
    LS_1 = LX + \frac{1}{4\cdot 3}LX^4 \equiv -X^2 + \frac{1}{4\cdot 3}(4\cdot 3 X^2 - X^5) \equiv - \frac{1}{4\cdot 3}X^5 \mbox{ mod } \partial,
\]
we have $\Phi(S_1)=-\frac{1}{4\cdot 3}X^5$ and so
\[
    S_2 = (1-\Phi'^{-1}\circ\Phi)s_1 = \left(X + \frac{1}{4\cdot 3}X^4\right) - \Phi'^{-1}\left(-\frac{1}{4\cdot 3}X^5\right) = X + \frac{1}{4\cdot 3}X^4 + \frac{1}{(4\cdot 3)(7\cdot 6)}X^7.
\]
Inductively, we obtain an $n-$th Newton\rq{s} approximation to be
\[
    S_n = \sum_{k=0}^n{\frac{3^k(\frac{2}{3})_k}{(3k+1)!}X^{3k+1}}.
\]
From this we obtain the another well-known linearly independent power series solution
\[
    Y_2(X) = \sum_{k=0}^\infty{\frac{3^k(\frac{2}{3})_k}{(3k+1)!}X^{3k+1}}.
\]

\end{enumerate}

\subsection{Distributional solutions of LDEs  as $\partial-$adic expansions and immediate Maps}\label{lde1} \label{Dis}

In this section, we equip $\mathcal{D}$ with the {\it unusual order} so that
$F^0\mathcal{D}=\mathbb{C}[\partial]$ (see Example \ref{Ex1}).

Similar to \S\ref{SS:lde}, let $L$ be an irreducible differential operator
of order $n$. $p(\partial)\in\mathbb{C}[\partial]$ is a polynomial so that the
order of $L-p(\partial)X^n$ is at most $n-1$. Let
$\Phi:\mathbb{C}[\partial]\to\mathbb{C}[\partial]$ be the remainder map of
$L/X$. Equip $\mathbb{C}[\partial]$ with the $\partial$-adic valuation $\nu$.
Such a map $\Phi$ becomes a map of ultrametric spaces.

Similar to the treatment in the \S\ref{SS:lde}, if $\nu$ denotes the $\partial$-adic
valuation, we infer from the estimate
\[
\nu(\Phi(f(\partial)))\geq\nu(f(\partial))-n\hspace{1cm}\mbox{for all }
f(\partial)\in\mathbb{C}[\partial]
\]
that $\Phi$ is continuous, and thus it is extended to the completion
$\overline{\mathbb{C}[\partial]}=\mathbb{C}[[\partial]]$ and we call the
extension $\Phi$ again.

Now we want to show that $\Phi$ is surjective by seeing that it is immediate.

Assume further that $p(0)\neq 0$. Define
$\Phi':\mathbb{C}[\partial]\to\mathbb{C}[\partial]$ by
\[
\Phi'(f(\partial))=p(0)f^{(n)}(\partial)\hspace{1cm}\mbox{for all }
f\in\mathbb{C}[\partial].
\]
Then $\Phi'$ is obviously immediate. Similar to the last section, we verify that
\[
\nu(\Phi(f(\partial))-\Phi'(f(\partial)))>\nu(\Phi(f(\partial)))\hspace{1cm}
\mbox{for all }f(\partial)\in\mathbb{C}[\partial].
\]
Therefore the conditions of Theorem \ref{T:prederivative} are all satisfied. Thus,
$\Phi$ is immediate (and hence surjective). There exists
$S\in\mathbb{C}[[\partial]]$ such that $\Phi(S)=0$. Such $S$ yields
\[
\overline{\mathcal{D}/\mathcal{D}L}\stackrel{\times S}{\longrightarrow}
\overline{\mathcal{D}/\mathcal{D}\partial}
\]
which is a ``distributional (supported at $0$) solution" of $L$.

\begin{remark}
Everything in this subsection can be obtained in the same way as that in the
previous subsection by interchanging the roles of $X$ and $\partial$. In fact,
the map from $\mathcal{D}$ to $\mathcal{D}$ defined by
\[
X\mapsto\partial,\hspace{1cm}\partial\mapsto-X
\]
is an automorphism of $\mathbb{C}$-algebra, which is called the Fourier transform.
Therefore, the results in this subsection can be obtained from the previous one
by applying the Fourier transform.
\end{remark}

\section{Regular Singularity and Immediate Maps}
In the last two examples, we studied the expansions of solutions of differential
equations at an ``ordinary point''. The next target is to study that at a \textit{regular singularity} that leads to Frobenius solutions of differential equations and difference equations as applications.
\bigskip

\subsection{Frobenius Series as an $X-$adic expansion and immediate Maps} \label{SS:Reg}

\begin{definition}[(\textbf{Order})] Let $\mathcal{D}_1=\mathbb{C}[X\partial,X]$ be a subalgebra of $\mathcal{D}$. Note that $[X\partial,X]=X$. We equip $\mathcal{D}_1$ with an order by making $\mathbb{C}[X]$ the set of elements of order zero.
\end{definition}

Thus axioms (O1-O3) in Definition \ref{D:order} are obviously satisfied.
\bigskip

\begin{definition}[\textbf{(Regular singularity)}] Let the order of $L\in \mathcal{D}_1$ be $n$ so that it is expressed as
\begin{equation} \label{ode}
L=p_n(X)(X\partial)^n+p_{n-1}(X)(X\partial)^{n-1}+...+p_0(X)
\end{equation}
for some $p_j\in\mathbb{C}[X]$, and $p_n(X)\neq 0$. We say that $L$ has a regular singularity at $0$ if $p_n(0)\neq 0$.
\end{definition}
\bigskip

\begin{definition}[(\textbf{Indicial polynomial})] We define
\begin{equation}\label{E:IP}
		i(\lambda)=p_n(0)\lambda^n+p_{n-1}(0)\lambda^{n-1}+\cdots +p_1(0)\lambda+
	p_0(0)
	\end{equation}
to be the \textit{indicial polynomial} of $L$ \eqref{ode} at $X=0$.
\end{definition}
\bigskip

Let $\lambda\in\mathbb{C}$ be a zero of the indicial polynomial \eqref{E:IP}
	\begin{equation}
		i(\lambda)=p_n(0)\lambda^n+p_{n-1}(0)\lambda^{n-1}+\cdots +p_1(0)\lambda+
	p_0(0)=0
	\end{equation}such that $i(\lambda+m)\neq 0$ for all $m\in\mathbb{N}$.

We intend to study the remainder map
$\Phi:\mathbb{C}[X]\to\mathbb{C}[X]$ of $L\div (X\partial-\lambda)=L/ (X\partial-\lambda)$. 

Note that the expression 
$X\partial-\lambda$ is ``monic'' in the sense that
\[
\mathbb{C}[X\partial,X]/\mathbb{C}[X\partial,X](X\partial-\lambda)\cong
\mathbb{C}[X]=F^0\mathbb{C}[X\partial,X]=F^0\mathcal{D}_1,
\] according to Definition \ref{D:monic}, so that it makes sense to speak of the remainder map above.
\bigskip

One endows $\mathbb{C}[X]$ with the $X$-adic valuation $\nu$. A moment of thought suggests that the estimate 
	\begin{equation}\label{E:valuation_xpartial}
\nu(\Phi(f(X)))\geq\nu(f(X))\hspace{1cm}\mbox{for all }f(X)\in\mathbb{C}[X],
	\end{equation} holds and from which one deduces immediately that the $\Phi$ is 
a continuous map.
\bigskip

Clearly we can extend $\Phi$ to the completion $\mathbb{C}[[X]]$ and we still call the extension by the notation $\Phi$.

\begin{lemma} Let $i(\lambda+m)\not=0$ for all $m\in\mathbb{N}$. Then
image$(\Phi)\subset X\mathbb{C}[[X]]$, so that $\Phi$ is regarded as a map
$\Phi:\mathbb{C}[[X]]\to X\mathbb{C}[[X]]$.
\end{lemma}

\begin{proof}
Since $LS=Q(X\partial -\lambda)+R$, so we have the remainder map
	\begin{equation}\label{E:remainder_xpartial}
		\begin{split}
		\Phi(S(X)) &=S(0) i(\lambda)\\
		&=S(0)\big(p_n(0)\lambda^n+p_{n-1}(0)\lambda^{n-1}+\cdots +p_1(0)\lambda+
	p_0(0)\big)=0\quad \mod X.
	\end{split}
	\end{equation}
This completes the proof.
\end{proof}
\bigskip

We now construct an appropriate \textit{tangential approximation map} for use in this section.
\bigskip

\begin{lemma}[(\textbf{Tangential map})] \label{L:tangent_map_xpartial}
Let $i(\lambda+m)\not=0,\ m\in\mathbb{N}$ and let $\Phi':\mathbb{C}[X]\to X\mathbb{C}[X]$ be the $\mathbb{C}$-linear map
	\[
		\Phi'(X^m) :=[p_n(0)(m+\lambda)^n+p_{n-1}(0)(m+\lambda)^{n-1}+...+p_0(0)]X^m
\hspace{1cm}\mbox{for all }m\in\mathbb{N}.
\]
Then for every $f(X)\in\mathbb{C}[X]$,
	\begin{equation}\label{E:rs_linear_map}
\nu(\Phi(f(X))-\Phi'(f(X)))>\nu(\Phi^\prime(f(X))).
	\end{equation}
\end{lemma}
\bigskip

\begin{proof}
For each $f(X)\in\mathbb{C}[X]$, write $f(X)=X^mg(X)$ for some $m\in\mathbb{N}$
and $g(0)\neq 0$ (so that $\nu(f(X))=m$.). Then
\begin{eqnarray*}
\Phi(f(X))&=&\Phi(X^mg(X))\\
          &=&g(0)[p_n(0)(m+\lambda)^n+p_{n-1}(0)(m+\lambda)^{n-1}+...+p_0(0)]X^m
             \;\;\;(\mbox{mod } X^{m+1})\\
          &=&\Phi'(f(X))\;\;\;(\mbox{mod } X^{m+1}).
\end{eqnarray*}Hence $\nu(\Phi(f(X))-\Phi'(f(X)))\geq m+1>m=\nu(f(X))=
\nu(\Phi(f(X)))=\nu(\Phi^\prime(f(X)))$ since we have equality in \eqref{E:valuation_xpartial} and 
$i(\lambda+m)\not=0$).
\end{proof}
\bigskip

We study the tangential map $\Phi^\prime$ of $\Phi$.

\begin{lemma}[(\textbf{Tangential map is immediate})]\label{L:Tangential map is immediate}
If $i(\lambda+m)\neq 0$ for all $m\in\mathbb{N}$, then the tangential map
$\Phi':\mathbb{C}[[X]]\to X\mathbb{C}[[X]]$ from the previous lemma is immediate.
\end{lemma}
\bigskip

\begin{proof}
For each $f(X)\in\mathbb{C}[[X]]$, we show that $Xf(X)=\sum_{k=1}^{\infty}f_kX^k$
is an attractor of $\Phi'$.
\bigskip

\begin{enumerate}
	\item[] \textbf{(AT1)} For each $G(X)\in\mathbb{C}[[X]]$ such that $\Phi'(G(X))\neq Xf(X)$, define
		\[
			F(X)=G(0)+ \sum_{k=1}^{\infty}\dfrac{f_k}{i(\lambda+k)}X^k.
		\]
Then the Lemma \eqref{L:tangent_map_xpartial} asserts that
		\[
			\nu(\Phi'(F(X))-Xf(X))=\infty>\nu(\Phi'(G(X))-Xf(X))
		\]
	\item[] \textbf{(AT2)} For each $H(X)\in B(F(X),\, G(X))$, we write
\[
\begin{array}{l}
G(X)=F(X)+X^mg(X)\\
H(X)=F(X)+X^kh(X)\hspace{1cm}\mbox{for some }k\geq m\mbox{ and }
g(X),h(X)\in\mathbb{C}[[X]],g(0)\neq 0.
\end{array}
\]
	\end{enumerate}
Then,
\[
\nu(\Phi'(H(X))-Xf(X))=k\geq m=\nu(\Phi'(G(X))-Xf(X)),
\]
so that $\Phi'(B(G(X),F(X)))\subset B(\Phi^\prime(G(X)),Xf(X))$.
\end{proof}
\bigskip

\begin{corollary} [(\textbf{Remainder map is immediate})] \label{C:Remainder map is immediate}
The remainder map $\Phi:\mathbb{C}[[X]]\to X\mathbb{C}[[X]]$ defined in  \eqref{E:remainder_xpartial} is immediate and hence surjective.
\end{corollary}
\bigskip

\begin{proof}
The completeness of $\mathbb{C}[[X]]$ together with the last two lemmas imply that
all conditions of Theorem \ref{T:prederivative} are satisfied. So $\Phi$ is
immediate.
\end{proof}
\bigskip

In particular, there exists non-zero $S\in\mathbb{C}[[X]]$ such that $\nu(S(X)-1)\geq 1$ and
$\Phi(S)=0$, which induces
\[
\overline{\mathbb{C}[X\partial,X]/\mathbb{C}[X\partial,X]L}
\stackrel{\times S}{\longrightarrow}
\overline{\mathbb{C}[X\partial,X]/\mathbb{C}[X\partial,X](X\partial-\lambda)}.
\]
This is a ``Frobenius series'' solution of the differential equation $L$.
\bigskip

The following theorem summarizes the discussion of this section. In addition, it gives the radius of convergence of the Frobenius series in parallel with the theorem for the power series counterpart for an ordinary point.
\bigskip

\begin{theorem} [(\textbf{Radius of convergence})]\label{T:radius_convergence}
Let $S=\sum_{k=0}^{\infty}a_kX^{k}$ be the $X-${adic} expansion of an element in the completion $\mathbb{C}[[X]]$ so that
	\[
\overline{\mathbb{C}[X\partial,X]/\mathbb{C}[X\partial,X]L}
\stackrel{\times S}{\longrightarrow}
\overline{\mathbb{C}[X\partial,X]/\mathbb{C}[X\partial,X](X\partial-\lambda)}
	\]is left $\mathcal{D}$-linear. Then
	\[
		\limsup_k|a_k|^{1/k}\leq 1/R
	\]where $R=\min\{|a|:p_n(a)=0\}$.
\end{theorem}
\bigskip

\subsection{Application to Bessel operators and immediate maps} \label{SS:bessel}

Suppose  that $-\nu\notin\mathbb{N}$.
We illustrate our theory by constructing a well-known series solution to Bessel equation. We emphasis that since our approach is completely general, it does not only apply to differential equations, but also to difference equations or $q-$difference equations. 

Since the Bessel operator 

	\begin{equation}\label{E:bessel}
		Ly=x^2y^{\prime\prime}+xy^{\prime}+(x^2-\nu^2)y
	\end{equation}possesses a regular singularity at $x=0$, we derive the well-known Bessel function of the first kind $J_\nu$ (Frobenius series) solution by computing the remainder maps and tangential maps.

Let us start with an initial approximate solution $S_0=1$ to compute for the remainder map $\Phi(1)$:
	\[
		\begin{split}
		L1& =X^2\partial^2+X\partial +(X^2-\nu^2)\\
			&=(X\partial)^2+(X^2-\nu^2)\\
			&=(X\partial+\nu)(X\partial-\nu) +X^2\\
			&=Q (X\partial-\nu)+\Phi(1)
		\end{split}
	\]where $Q=X\partial+\nu$ and $\Phi(1)=X^2$. Since 
	\[
		L= (X\partial)^2+(X^2-\nu^2)
	\]so that the tangential map \eqref{L:tangent_map_xpartial} and  the indicial polynomial \eqref{E:IP} becomes
	\[
		\Phi^\prime(X^m)=[(m+\nu)^2-\nu^2]=(m^2+2m\nu)X^m,\quad m\in\mathbb{N}.
	\]Hence 
	\[
		\Phi(X^2)=2^2(1+\nu)X^2.
	\]Since
	\[
		X^2={\Phi^{\prime}}^{-1}\circ \Phi(X^2)={\Phi^{\prime}}^{-1}(2^2(1+\nu)X^2)
		=2^2(1+\nu){\Phi^{\prime}}^{-1}(X^2),
	\]so
	\[
		{\Phi^{\prime}}^{-1}(X^2)=\frac{X^2}{2^2(1+\nu)}.
	\]
Next we compute the first improved root $S_1$:
	\[
		\begin{split}
			S_1 &=1-{\Phi^{\prime}}^{-1}\circ \Phi(S_0)\\
				&=1-{\Phi^{\prime}}^{-1}\circ \Phi(1)=1-\frac{X^2}{2^2(1+\nu)}.
		\end{split}
	\]It is straightforward to verify that the formula
	\begin{equation}\label{E:bessel_intermediate}
		LX^{2k}=[(X\partial+\nu)(X\partial-\nu) +X^2]X^{2k}
		=2^2k(k+\nu)X^{2k}+X^{2k+2},
	\end{equation}holds for every $k\in\mathbb{N}$. 
	
	We apply this formula to compute for the second remainder map:
		\[
			\begin{split}
				LS_1 &= [(X\partial+\nu)(X\partial-\nu) +X^2]S_1\\
					&=[(X\partial+\nu)(X\partial-\nu) +X^2]\big(1- \frac{X^2}{2^2(1+\nu)}\big)\\
					&\equiv-\frac{X^4}{2^2(1+\nu)}\quad \mod X\partial-\nu.
			\end{split}
		\]It is easy to check that
		\begin{equation}\label{E:inverse_tangential}
			{\Phi^\prime}^{-1} (X^m)=\frac{X^m}{m^2+2m\nu}, \quad m\in \mathbb{N}
		\end{equation}holds. Hence the second improved root $s_2$ is given by
		\[
			\begin{split}
			S_2 &= (1-{\Phi^\prime}^{-1}\circ \Phi)(S_1)\\
				&= \big(1-\frac{X^2}{2^2(1+\nu)}\big)-({\Phi^\prime}^{-1}\circ \Phi)
	\big(1-\frac{X^2}{2^2(1+\nu)}\big)\\
				&=  \big(1-\frac{X^2}{2^2(1+\nu)}\big)-{\Phi^\prime}^{-1} \big(-\frac{X^4}{2^2(1+\nu)}\big) \quad \mod X\partial-\nu\\
				&=1-\frac{X^2}{2^2(1+\nu)}+\frac{X^2}{2^2(1+\nu)}\frac{X^4}{16+8\nu} \quad \mod X\partial-\nu\\
				&=1-\frac{(X/2)^2}{(1+\nu)\cdot 1!}+\frac{(X/2)^{2\cdot 2}}{(1+\nu)(2+\nu)\cdot 2!},
 \quad \mod X\partial-\nu.
			\end{split}
		\]
		We proceed to the inductive step to assume that we already have derived the $n-$th improved root
		\[
			S_n=\sum_{k=0}^n\frac{(-1)^k (X/2)^{2k}}{(1+\nu)\cdots (k+\nu)\cdot k!}.
		\]Let us now compute the $(n+1)-$th improved root
		\[
			S_{n+1}=[1-{\Phi^\prime}^{-1}\circ \Phi](S_n).
		\]To do so, we first compute the corresponding remainder map $\Phi(S_n)$:
		\[
			\begin{split}
				LS_n &=[(X\partial+\nu)(X\partial-\nu) +X^2]S_n\\
					&=[(X\partial+\nu)(X\partial-\nu) +X^2] \big(\sum_{k=0}^n\frac{(-1)^k (X/2)^{2k}}{(1+\nu)\cdots (k+\nu)\cdot k!}\big)\\
					&=\sum_{k=0}^n\frac{(-1)^k [(X\partial+\nu)(X\partial-\nu)+x^2]X^{2k}}{2^{2k}(1+\nu)\cdots (k+\nu)\cdot k!}
\\
					&=\sum_{k=0}^n\frac{(-1)^k (2^2k(k+\nu)X^{2k})}{2^{2k}(1+\nu)\cdots (k+\nu)\cdot k!}+\sum_{k=0}^n\frac{(-1)^k X^{2k+2}}{2^{2k}(1+\nu)\cdots (k+\nu)\cdot k!}\\
					&=\sum_{k=1}^n\frac{(-1)^k X^{2k}}{2^{2k-2}(1+\nu)\cdots (k-1+\nu)\cdot (k-1)!}+\sum_{k=1}^{n+1}\frac{(-1)^{k-1} X^{2k}}{2^{2k-2}(1+\nu)\cdots (k-1+\nu)\cdot (k-1)!}\\
					&=\frac{(-1)^n X^{2n+2}}{2^{2n}(1+\nu)\cdots (\nu+n)n!},
					\qquad \mod x\partial-\nu\\
					&=\Phi(S_n).
			\end{split}
		\]
We deduce immediately that
	\[
		\begin{split}
			S_{n+1} &= [1-{\Phi^\prime}^{-1}\circ \Phi](S_n)\\
			&=S_n-{\Phi^\prime}^{-1} \Phi(S_n)\\
			&=\sum_{k=0}^n\frac{(-1)^k (X/2)^{2k}}{(1+\nu)\cdots (k+\nu)\cdot k!}-
			{\Phi^\prime}^{-1}\Big(\frac{(-1)^n X^{2n+2}}{2^{2n}(1+\nu)\cdots (\nu+n)n!}\Big)
			\qquad  \mod x\partial-\nu\\
			&=\sum_{k=0}^n\frac{(-1)^k (X/2)^{2k}}{(1+\nu)\cdots (k+\nu)\cdot k!}-
			\frac{(-1)^n X^{2n+2}}{2^{2n}(1+\nu)\cdots (\nu+n)n!}\frac{1}{(n+1)(n+1+\nu)}\\
			&=\sum_{k=0}^{n+1}\frac{(-1)^k (X/2)^{2k}}{(1+\nu)\cdots (k+\nu)\cdot k!},
		\end{split}
	\]as desired. Since $-\nu\notin\mathbb{N}$, so it follows from Lemma \ref{L:Tangential map is immediate}, Corollary \ref{C:Remainder map is immediate} that the $\Phi^\prime$ and $\Phi$ are respectively immediate. We deduce that the expression
	\begin{equation}\label{E:Bessel_eqn}
	\sum_{k=0}^{\infty}\frac{(-1)^k (X/2)^{2k}}{(1+\nu)\cdots (k+\nu)\cdot k!}	
	\end{equation}in $\mathbb{C}[[x]]$
  	provides an $X-$adic expansion as a solution map in
\[
	\overline{\mathbb{C}[X\partial,X]/\mathbb{C}[X\partial,X]L}
	\stackrel{\times S}{\longrightarrow}
	\overline{\mathbb{C}[X\partial,X]/\mathbb{C}[X\partial,X](X\partial-\nu)}
\]
where
	\[
		L= (X\partial)^2+(X^2-\nu^2).
	\]			
Theorem \ref{T:radius_convergence} immediate suggests that the series \eqref{E:Bessel_eqn} has infinite radius of convergence $R=+\infty$ when it is interpreted in the complex plane $\mathbb{C}$. 

Similar consideration of another indicial root $-\nu$ when $\nu\notin\mathbb{Z}$ would give rise to a solution of $(X\partial)^2+X^2-\nu^2$ in the $X-adic$ completion of $\mathcal{D}/\mathcal{D}(X\partial+\nu)$. This would give a second linearly independent Frobenius solution to the Bessel operator. 
\bigskip

\begin{remark} We note that if $\nu\in\mathbb{Z}$, then one needs to construct $\mathcal{D}-$linear maps for
	\[
	\overline{\mathbb{C}[X\partial,X]/\mathbb{C}[X\partial,X]L}
	\stackrel{\times S}{\longrightarrow}
	\overline{\mathbb{C}[X\partial,X]/\mathbb{C}[X\partial,X](X\partial-\nu)^2}
	\] by studying its remainder map, which is reserved for a future project.
\end{remark}
\bigskip

\subsection{Application to algebraic Huen operators and immediate maps}\label{SS:pseries}
To study a power series solution $f=\sum_{m=0}^\infty b_mx^m$ of Heun equation $\mathscr{H}f=0$ at the regular singular point $x=0$, instead of considering solutions in the space $\mathcal{D}/\mathcal{D}(X\partial)$,  as in the study in \S\ref{SS:Reg}, we may also consider solutions in the space $\mathcal{D}/\mathcal{D}\partial$ as below.

\begin{proposition}\label{P:remainder}
For each a polynomial $S\in\mathbb{C}[X]$ of degree $n$, there exist polynomials $q_S,r_S\in\mathbb{C}[X]$ with
$\deg q_S \leq n+2$, $\deg r_S \leq n+1$ such
that
\[\mathscr{H} S=[S(X-1)(X-a)\partial+q_S]\partial+r_S.
	\]
\end{proposition}

\begin{proof}[Sketch of Proof]
By linearity, we suffice to show the statement is true when $S=X^n$ for any $n=0,1,2,\cdots$.
We prove it by induction. First, we can see that
\[\mathscr{H}1=\{X(X-1)(X-a)\partial+(\gamma+\delta+\epsilon)X^2-[\gamma+\epsilon+(\gamma+\delta)a]X+\gamma a\}\partial+\alpha\beta X.\]
Then by induction and the property that $[\partial,X]=1$, for any $m=0,1,\cdots,$ 
\begin{eqnarray*}
\mathscr{H}X^{m+1}&=&\{X^{m+2}(X-1)(X-a)\partial+A_mX^{m+3}+B_mX^{m+2}+C_mX^{m+1}\}\partial\\
                      & &\  +D_mX^{m+2}+E_mX^{m+1}+F_mX^{m},
\end{eqnarray*}
where
	\[
		\begin{split}
			A_{m} & =\alpha+\beta+2m+3;\\
			B_{m} &=-[\gamma+\epsilon+2m+2+(\gamma+\delta+2m+2)a];\\
			C_m &=(\gamma+2m+2)a;\\
			D_{m} &=(\alpha+m+1)(\beta+m+1);\\
			E_{m} &=-(m+1)[\gamma+\epsilon+m+(\gamma+\delta+m)a];\\
			F_{m} &=(\gamma+m)(m+1)a.
	\end{split}
	\]
\end{proof}
\bigskip

By the above proposition, it makes sense to define
\begin{definition}\label{E:heun_remainder}
Let the remainder map $\Psi:\mathbb{C}[X]\rightarrow\mathbb{C}[X]$ be defined by $\Psi(S)=r_S$ for $s\in\mathbb{C}[X]$.
\end{definition}
\bigskip

\begin{corollary}
Let $\nu:\mathbb{C}[X]\to\mathbb{N}\cup\{\infty\}$ be the $X$-adic valuation. Define the tangential map be $\mathbb{C}$-linear map
$\Psi':\mathbb{C}[X]\to\mathbb{C}[X]$ by
\[
\Psi'(X^m)=m(\gamma+m-1)aX^{m-1}  \hspace{1cm}\mbox{for all }m\in\mathbb{N}.
\]
Then for each $f(X)\in\mathbb{C}[X]$,
\[
\nu(\Psi(f(X))-\Psi'(f(X)))>\nu(\Psi(f(X)))=\nu(\Psi^\prime (f(X))).
\]
\end{corollary}
\bigskip

\begin{corollary} The remainder map 
$\Psi:\mathbb{C}[[X]]\to\mathbb{C}[[X]]$ is immediate and hence surjective.
\end{corollary}
\bigskip

\begin{proof}
Observe that $\Psi'$ is immediate. The last corollary says that the conditions in Theorem \ref{T:prederivative} are satisfied.
Thus, $\Psi$ is immediate and hence surjective.
\end{proof}
\bigskip

Then it follows that there exists a non-zero $s_0\in\mathbb{C}[[X]]$ such that $\Psi(S_0)=0$. In other words,
\[
\mathscr{H}S_0=[S_0(X-1)(X-a)\partial+q_{S_0}]\partial\hspace{1cm}\mbox{for some }q.
\]

\begin{remark}
The power series
$S_0=\sum_{m=0}^\infty c_mX^m$ satisfies the Heun equation \[\mathscr{H}S_0=Q_{S_0}\partial(1)=0.\]
\end{remark}

\begin{remark}
Instead of $\mathscr{H}$, the above argument still holds for $\mathscr{H}-q$ for any $q\in\mathbb{C}$.
\end{remark}		

\medskip


\subsection{${}_1F_1$-Series as $\partial-$adic expansions of solutions of doubly confluent Heun equations} \label{SS:DHE}

As an application of the method developed so far, we will study series solutions to the well-known \textit{doubly
confluent Heun equation} \cite{Ronveaux} from the viewpoint taken in this article. We first fix our notation.

\begin{definition}
Fix $a$, $b$, $c$, $q\in\mathbb{C}$.  We denote the well-known {\it confluent hypergeometric operator} by
\begin{equation} \label{E:CH2}
K=X\partial^2+(b-X)\partial-a,
\end{equation}
where $-a\not\in \mathbb{N}$ so that $K$ is irreducible. We also denote 
{\it doubly confluent Heun operator} by 
	\begin{equation} \label{E:DHE}
	D=X^2\partial^2+(-X^2+bX+c)\partial-aX+q=XK+c\partial+q.
	\end{equation}
\end{definition}

Thus our target is to solve the doubly confluent Heun equation in the left
$\mathcal{D}$-module $\mathcal{D}/\mathcal{D}K$. 
For this setup, we define an unconventional kind of \textit{order} in $\mathcal{D}$ so that 
 $F^0\mathcal{D}=\mathbb{C}[\partial]$ as in Definition \ref{D:order}. In other words the  \textit{order} of an element in $\mathcal{D}$ is the degree of $X$ involves there, as discussed in Example \ref{Ex1}. The confluent hypergeometric operator $K$ is therefore a differential operator of order $1$ under this interpretation.

It is clear that $K$ is not monic under the Definition \ref{D:monic}, i.e., the composition
	\[
		F^0\mathcal{D}\longrightarrow\mathcal{D} \longrightarrow\mathcal{D}/\mathcal{D}K
	\]is clearly \textit{not} a $\mathbb{C}$-vector space isomorphism. However, the ``leading coefficient" of $K$ \eqref{E:CH2}, being $\partial^2-\partial$, happens to be the same as that of $D$ \eqref{E:DHE}. Thus we expect the remainder map
$\Phi:\mathbb{C}[\partial]\to\mathbb{C}[\partial]$ of $D/K$ to be well-defined.
In fact, it is a consequence of the following commutator relation.
\bigskip

\begin{lemma}\label{L:K-partial}
Let $n\in\mathbb{N}\cup\{0\}$, we have
	\begin{enumerate}
		\item
		\[
			[\partial^n,\, K]=n(\partial^{n+1}-\partial^n).
	\]
		\item 
	\[
		D\partial^n=(X\partial^n-n\partial^{n-1})K+c\partial^{n+1}+(nb+n(n-1)+q)\partial^n
-n(a+n-1)\partial^{n-1}.
\]
	\end{enumerate}
\end{lemma}
\bigskip

\begin{proof}Direct computation.
\end{proof}
\bigskip

Now we endow $\mathbb{C}[\partial]$ with the
$\partial$-adic valuation $\nu$. The last lemma above yields the following estimate.
\bigskip

\begin{lemma}
For each $f(\partial)\in\mathbb{C}[\partial]$, the valuation of the remainder map $\Phi$ satisfies the inequality
\[
\nu(\Phi(f(\partial)))\geq\nu(f(\partial))-1.
\]
\end{lemma}
\bigskip

In particular, $\Phi$ is continuous and thus it can be extended to the $\partial$-adic
completion $\mathbb{C}[[\partial]]$ and we call the extension $\Phi$ again. Similar to
the previous examples, the next step is to obtain an approximation of $\Phi$. Indeed the following is a direct consequence of Lemma \ref{L:K-partial}.
\bigskip

\begin{lemma}
Define a $\mathbb{C}$-linear tangential map $\Phi':\mathbb{C}[\partial]\to\mathbb{C}[\partial]$
so that for each $n\in\mathbb{N}$
\[
\Phi'(\partial^n)=-n(a+n-1)\partial^{n-1},
\]
where $-a\not\in \mathbb{N}$.
Then 
	\begin{enumerate}
		\item $\Phi^\prime$ is immediate
		\item  for each $f(\partial)\in\partial \mathbb{C}[\partial]$,
\[
\nu(\Phi(f(\partial))-\Phi'(f(\partial)))>\nu(\Phi(f(\partial)))=\nu(\Phi^\prime (f(\partial))).
\]
	\end{enumerate}
\end{lemma}
\bigskip

\begin{proof} Direct verification.
\end{proof}
\bigskip

The following theorem follows similarly.
\bigskip

\begin{corollary} Let $-a\notin\mathbb{N}$. Then $\Phi:\mathbb{C}[[\partial]]\to \mathbb{C}[[\partial]]$
is immediate and hence surjective.
\end{corollary}
\begin{proof}
Since $\Phi'$ is immediate if $-a\notin\mathbb{N}$ (see Definition \ref{D:immediate}), so $\Phi$ is immediate and hence surjective also by 
applying the last lemma and Theorem \ref{T:prederivative}.
\end{proof}
\bigskip

Hence there exists a non-zero $S\in\mathbb{C}[[\partial]]$ such that $\Phi(s)=0$.
Such $s$ yields a left $\mathcal{D}$-linear map
\[
\overline{\mathcal{D}/\mathcal{D}D}\stackrel{\times S}{\longrightarrow}
\overline{\mathcal{D}/\mathcal{D}K},
\]
which is an expansion of solutions of the double confluent operator $D$ as a series of $\partial^k {}_1F_1(a,b;\, x)$,
although the common practice is to write as a series of ${}_1F_1(a+k,b+k;\ x)$.

\begin{theorem} 
Let $S=\sum_{k=0}^{\infty}a_k \partial^{k}$ be the $\partial-${adic} expansion of an element in the completion $\mathbb{C}[[\partial]]$ so that
	\[
\overline{\mathbb{C}[X, \partial]/\mathbb{C}[X,\, \partial]D}
\stackrel{\times S}{\longrightarrow}
\overline{\mathbb{C}[X,\, \partial]/\mathbb{C}[X,\,\partial]K}.
	\]is left $\mathcal{D}$-linear. Then
	\[
		\limsup_k|a_k|^{1/k}=0.
	\]
\end{theorem}


\subsection{Parabolic cylindrical series as $(\partial\pm X)-$adic expansions of solutions of the biconfluent Heun equations} \label{SS:BHE}

\begin{definition}
Let the {\it creation} and {\it annihilation} operators be defined by
\[
A^\dagger=\partial-X\quad  \mbox{ and }\quad  A=\partial+X
\]
respectively. The {\it simple harmonic operator} is defined by 
	\[
		H=A^\dagger A-1=AA^\dagger +1.
	\]
\end{definition}
\bigskip

We obviously have
\begin{lemma}
\[
[A^\dagger,A]=2.
\]
\end{lemma}
\bigskip

\begin{definition}
Fix $\alpha$, $\beta$ and $\gamma\in\mathbb{C}$. The (normalized) {\it biconfluent
Heun operator} is
\begin{equation} \label{biconfluent}
	B=(A-A^\dagger)^2(H+\alpha)+\beta(A-A^\dagger)+\gamma.
\end{equation}
We will suppress $\alpha$, $\beta$ and $\gamma$ in the symbol $B$ for simplicity. 
\end{definition} 
\bigskip

The normalized form of the \textit{biconfluent Heun equation}  assumes the form (see e.g., \cite[p. 236]{DL}, \cite{Ronveaux}, \cite[(5.2)]{Chiang_Yu})
	\begin{equation}\label{E:BHE}
		y^{\prime\prime}-\big(x^2+\beta x+\frac{\beta^2}{4}-\gamma+\frac{\delta}{2x}+\frac{\alpha^2-1}{4x^2}\big)y=0.
	\end{equation}
We notice that the parameter $\beta$ in \eqref{E:BHE} can be removed without affecting the study of its solutions. Thus without loss of generality, we may study the \eqref{biconfluent} which contains three parameters $\alpha,\, \beta,\, \gamma$ (after relabeling)  instead of the four parameters in \eqref{E:BHE}.
\bigskip

\begin{definition} Let $F^0\mathcal{D}=\mathbb{C}[A^\dagger]$ be the ring of zeroth order operators.  That is, for each differential operator in $\mathbb{C}[A^\dagger,A]/<[A^\dagger,A]-2>$, its unconventional
{\it order} is the highest power of $A$ involved.
\end{definition}
\bigskip

Note that the biconfluent Heun operator can be rewritten as
\[
\begin{array}{rl}
A^{\dagger 2}B=&A^{\dagger 2}(A^2-A^\dagger A-AA^\dagger+A^{\dagger 2})(A^\dagger A-1+\alpha)
 +A^{\dagger 2}\beta(A-A^\dagger)+\gamma A^{\dagger 2}\\
 =&(A^\dagger A)^3+[-2A^{\dagger 2}+\alpha+1](A^\dagger A)^2\\
 &+[A^{\dagger 4}+(4-2\alpha)A^{\dagger 2}+\beta A^\dagger+2(\alpha-1)](A^\dagger A)\\
 &+(\alpha-1)A^{\dagger 4}-\beta A^{\dagger 3}+(2\alpha-2+\gamma)A^{\dagger 2}.
\end{array}
\]
Note that the expression $A^\dagger A$ plays the role of $X\partial$ in the study of regular singularity earlier. 
This operator has a regular singular point at $0$ whose indicial polynomial (see \eqref{E:IP}) is
\[
\lambda^3+(\alpha+1)\lambda^2+2(\alpha-1)\lambda=\lambda(\lambda+2)(\lambda+\alpha-1).
\]
Obviously, $-\alpha+1$ is a zero of this polynomial. By our previous study of
operators with regular singular point at $0$, we study the remainder map
of $A^{\dagger 2}B/(H+\alpha)$.
However, a simple verification implies that
\[
[A^{\dagger 2},B]\in\mathcal{D}(H+\alpha)+\mathbb{C}[A^\dagger].
\]
It is equivalent to study the remainder map
$\Phi:\mathbb{C}[A^\dagger]\to\mathbb{C}[A^\dagger]$ of
$BA^{\dagger 2}/(H+\alpha)$. By the result about regular singular points as discussed in \S\ref{SS:Reg},
 we deduce 
\bigskip

\begin{theorem} There exists a non-zero $S=\sum_{k=2}^{\infty}a_k {A^\dagger}^k\in {A^\dagger}^2\mathbb{C}[[A^\dagger]]$  
 which induces a left
$\mathcal{D}$-linear map
\[
\overline{\mathcal{D}/\mathcal{D}B}\stackrel{\times S}{\longrightarrow}
\overline{\mathcal{D}/\mathcal{D}(H+\alpha)}.
\]Moreover, 
\[
		\limsup_k|a_k|^{1/k}=0.
	\]
\end{theorem}
\bigskip

Swapping $A^\dagger$ by $A$, we have 

\[
\begin{array}{rl}
A^2B=&A^{2}(A^2-A^\dagger A-AA^\dagger+A^{\dagger 2})(AA^\dagger +\alpha+1)
 +A^{2}\beta(A-A^\dagger)+\gamma A^{2}\\
 =&(A^\dagger A)^3+[-2A^{\dagger 2}+\alpha+1](A^\dagger A)^2\\
 &+[A^{\dagger 4}+(4-2\alpha)A^{\dagger 2}+\beta A^\dagger+2(\alpha-1)](A^\dagger A)\\
 &+(\alpha-1)A^{\dagger 4}-\beta A^{\dagger 3}+(2\alpha-2+\gamma)A^{\dagger 2}.
\end{array}
\]

If $D_{\alpha}$ is the parabolic cylindrical function (which is defined on the
upper half plane) solving $H+\alpha$, then $SD_{\alpha}$ is a formal solution of
the biconfluent Heun equation.

If $S=\sum_{m=2}^{\infty}a_m(A^\dagger)^m$, it is interesting to see if the formal
series $SD_{\alpha}$ converges in the classical sense.

First of all, $0$ is the only singular point of $B$, so
\[
\limsup_m|a_m|^{1/m}=0.
\]
Pick any point $a$ in the upper half plane and $R<$Im $a$. Let
$M=\max\{D_{\alpha}(z):|z-a|=R\}$. Then if $|z-a|<R/2$,
\[
|(A^\dagger)^m D_{\alpha}(z)|\leq Mm!(|a|+R+\dfrac{2}{R})^m.
\]
But
\[
\limsup_m\left|Mm!a_m(|a|+R+\dfrac{2}{R})^m\right|^{1/m}=0.
\]
Consequently, the formal $SD_{\alpha}$ converges uniformly on every compact subset
of the upper half plane and so it defines an analytic function in the upper half
plane in the classical sense.

\section{Table of Main Computational Results of Classical Functions}

We have already worked out the $n-$th Newton-Raphson approximation $S_n$ for the zeros of the remainder maps $\Phi$ for the Airy operator and Bessel operator in \S\ref{SS:airy} and \S\ref{SS:bessel} above respectively. In addition, we have further worked out these Newton-Raphson approximation $S_n$ for the zeros of the remainder maps $\Phi$ for the Hermite operator, confluent hypergeometric operator and the Gauss hypergeometric operator. Moreover, we have also worked out some of these approximations with different initial trial approximation $S_0$ that lead to \textit{linearly independent multiplication maps} for the Airy and confluent hypergeometric operators. We skip the tedious computation that lead to the multiplication maps and simply list the results in the following tables for reference. Needless to say that these maps apply not only to the familiar classical special functions that satisfy differential equations in the complex plane, but also other interpolation series type maps that satisfy difference equations for difference operators $\partial$ that conform to $[\partial,\, X]=1$ and where the operator $K$ and maximal ideal $\mathcal{I}$ can be constructed.

However, the construction of the remainder map $\Phi$ for the Heun operator $\mathscr{H}$ is not as explicit as those in the table below for the classical special functions. This is because we could only ensure that the valuation $\nu$  to satisfy the equality
	\[
		n+1\le \nu (\Phi(S_n))\le 2n+2.
	\]
That is, the ``effective portion" of the multiplication map $S_n$ which has a degree much higher than $n$, for the Heun operator. We hope to return to this issue in the future.

\begin{landscape}
\bigskip

\subsection{Table of operator $L$, first order operator $K$ and maximal ideal $I$  of examples}

\bigskip

\begin{center}
\begin{tabular}{|l|l|l|l|}
  \hline
  {\bf Section} & {\bf Linear Differential Operator $L$} &
  {\bf First-Order Operator $K$} & {\bf Maximal Ideal $I$}\\
  \hline
  \S\ref{SS:lde} & Usual $n$-th order operator  & $\partial$ & $X\mathbb{C}[X]$\\
  \hline
  \S\ref{Dis} & Usual $n$-th order operator & $X$ & $\partial\mathbb{C}[\partial]$\\
  \hline
  \S\ref{SS:Reg} & see (\ref{ode}) & $X\partial-\lambda$ ($\lambda$ is an indicial root) & $X\mathbb{C}[X]$\\
  \hline
  \S\ref{SS:DHE} & $X^2\partial^2+(-X^2+bX+c)\partial-aX+q$ &  $X\partial^2+(b-X)\partial -a$
 & $\partial\mathbb{C}[\partial]$\\
  \hline
  \S\ref{SS:BHE} & $(A-A^\dagger)^2(H+\alpha)+\beta(A-A^\dagger)+\gamma$ & $H+\alpha=A^\dagger A+\alpha-1$& $A^\dagger\mathbb{C}[A^\dagger]$ or $A\mathbb{C}[A]$\\
  \hline
\end{tabular}
\end{center}
\bigskip

\subsection{Table of remainders maps, tangential maps of well-known functions}

\bigskip

\begin{center}
\begin{tabular}{|l|l|l|l|}
  \hline
  {\bf Operators $L$} 
  	& {\bf Remainder Maps $\Phi$} 
	& {\bf  {\small $({\Phi^{-1}})^\prime (X^m)$}} 
	& {\small \bf $S_0$; \ Newton\rq{s} approximation $S_n$}
\\
  \hline
  	{\small $\partial^2-X$}  
 	& {\small $m(m-1)X^{m-2}-X^{m+1}$}  
	& {\small $\frac{X^{m+2}}{(m+2)(m+1)}$}  
  	& {\small $X$; $\sum_{k=0}^n{\frac{3^k(\frac{2}{3})_k}{(3k+1)!}X^{3k+1}}$}\\
 \hline
  	{\small $\partial^2-X$}
	 &{\small $m(m-1)X^{m-2}-X^{m+1}$}
	&{\small $\frac{X^{m+2}}{(m+2)(m+1)}$ }
	&{\small $1$; $\sum_{k=0}^n{\frac{3^k(\frac{1}{3})_k}{(3k)!}X^{3k}}$}\\
  \hline
	{\small $\partial^2-2X\partial+2\lambda$}
	 &{\small $m(m-1)X^{m-2}+2(\lambda-m)X^m$}
	 &{\small $\frac{X^{m+2}}{(m+2)(m+1)}$}
	 &{\small $1$; $\sum_{k=0}^n\frac{(-1)^k2^k (\lambda)_{2k} }{(2k)!}X^{2k}$}\\
 \hline
	{\small $\partial^2-2X\partial+2\lambda$}
	&{\small $m(m-1)X^{m-2}+2(\lambda-m)X^m$}
	&{\small $\frac{X^{m+2}}{(m+2)(m+1)}$}
	& {\small $X$; $\sum_{k=0}^n\frac{(-1)^k2^k(\lambda-1)_{2k+1}}{(2k)!}X^{2k+1}$}\\
  \hline
	{\small $X(X\partial^2+(\alpha_1-x)\partial+n)$}
	&{\small $(n-m)X^{m+1}+m(\alpha+m)X^m$}  
	&{\small $ \frac{X^m}{m(\alpha+m)}$}
	& {\small $1$; $\sum_{k=0}^n \frac{(-a)_k}{(\alpha+1)_kk!}X^k$}\\
  \hline
	 {\small $X(X(1-X)\partial^2+[X-(a+b+1)X]\partial-ab)$}
	& {{\small $(c-1+m)mX^m -(a+m)(b+m)X^{m+1}$}}
	&{\small $\frac{X^m}{(c-1+m)m}$}
	& {\small $1$;\ $\sum_{k=0}^n \frac{(a)_k(b)_k}{(c)_k k!} X^k$} \\
  \hline
       {\small $X\mathscr{H}$}
       & {\small $D_{m-1}X^{m+2}+E_{m-1}X^{m+1}+F_{m-1}X^m$} 
       & {\small $ \frac{X^m}{F_{m-1}}$}
       & {\small $\nu(\Phi(X^n))=n+1$}\\
  \hline
\end{tabular}
\end{center}
\end{landscape}


\part{Eigenvalue problems of Remainder Maps and Applications}\label{P:2}

\section{Hypergeometric and Heun Differential Operators} \label{HE}

In the last section, we found factorizations of operators (in certain completions)
by dividing operators by a first order operator. The remainder is a zeroth order
operator and there is a way to judge if this remainder is small.

The same procedure does not appear to work if we divide operators by a second
order operator, as the remainder would be of first order. Therefore, the usual $I$-adic
ultrametric approach is not applicable to such a space of first order operators.

However, in case the dividend is of certain special types,
we can still guarantee the remainders to
have order zero. Then, the method of successive approximate root of an immediate
map is again applicable. The study of these special cases turn up to be fruitful
and it leads to a partial solution of a conjecture by Takemura for factorization of generalized hypergeometric operators. Moreover, our theory also extends the study to factorization of generalized confluent hypergeometric operators. The rest of
the article is devoted to such consideration.

Let the Heun operator $ \mathscr{H}$ be
  \begin{equation}\label{E:heun}
	\mathscr{H}(a;\alpha,\beta,\gamma,\delta):=
     X(X-1)(X-a)\partial^2+[\gamma(X-1)(X-a)+\epsilon X(X-1)+\delta X(X-a)]\partial
     + \alpha\beta X,
  \end{equation} with $\alpha+\beta-\gamma-\delta-\epsilon+1=0$ with the corresponding Riemann $P$-scheme \[
		P_{\mathbb{CP}^1}
		\begin{Bmatrix}
\ 0\ &\ 1\ &\ a\ &\ \infty\ &\\
0&0&0&\alpha&x \\
1-\gamma&1-\delta&1-\epsilon&\beta&
	\end{Bmatrix}.
	\]
We will study the solution of Heun equation by the series of hypergeometric functions as follows:
\[
	\sum_{m=0}^{\infty}C_m \cdot	P_{\mathbb{CP}^1}
		\begin{Bmatrix}
\ 0\ &\ 1\ &\ \infty\ &\\
0&0&\alpha&x \\
1-\gamma + m&1-\delta-\epsilon - m&\beta&
	\end{Bmatrix}
	\] and
\[
	\sum_{m=0}^{\infty}\hat{C}_m \cdot	P_{\mathbb{CP}^1}
		\begin{Bmatrix}
\ 0\ &\ 1\ &\ \infty\ &\\
0&0&\alpha&x \\
1-\gamma-\epsilon - m&1-\delta + m&\beta&
	\end{Bmatrix}.
	\]
	
Let the hypergeometric operators $\mathcal{H}$ and $\hat{\mathcal{H}}$(with $m=0$) be
	\begin{equation}\label{E:H1-Heun}
		 \mathcal{H}(\alpha,\beta,\gamma,\delta):=
     X(X-1)\partial^2+[\gamma(X-1)+(\delta+\epsilon)X]\partial
     + \alpha\beta
	\end{equation}
	 and
\begin{equation}\label{E:H2-Heun} 
\hat{\mathcal{H}}(\alpha,\beta,\gamma,\delta):=
     X(X-1)\partial^2+[(\gamma+\epsilon)(X-1)+\delta X]\partial
     + \alpha\beta,
\end{equation} with the Riemann $P$-schemes \[
		P_{\mathbb{CP}^1}
		\begin{Bmatrix}
\ 0\ &\ 1\ &\ \infty\ &\\
0&0&\alpha&x \\
1-\gamma&1-\delta-\epsilon&\beta&
	\end{Bmatrix}
\ \mbox{ and } \	
P_{\mathbb{CP}^1}
		\begin{Bmatrix}
\ 0\ &\ 1\ &\ \infty\ &\\
0&0&\alpha&x \\
1-\gamma-\epsilon & 1-\delta & \beta &
	\end{Bmatrix}.
	\]
\bigskip

	
\medskip

\subsection{Division of Heun operators by Hypergeometric operators}\label{hseries}	

In this section, we consider the ``annihilation'' operator $A$ and ``creation'' operator $A^\dagger$ which are defined by
	\begin{equation}\label{E:A1-Heun}
		A=X\partial+\gamma \quad \mbox{ and }\quad A^\dagger=(X-1)\partial+\delta+\epsilon-1=(X-1)\partial+\alpha+\beta-\gamma.
	\end{equation}

It is well-known (for example, see ($3'$) and ($6'$) in \cite[p. 93]{Poole}) that
	\[
A\left({}_2F_1\left(\begin{array}{c}\alpha,\ \beta\\ \gamma-m\end{array};\ x\right)\right)
=(\gamma-m-1){}_2F_1\left(\begin{array}{c}\alpha,\ \beta\\ \gamma-m-1\end{array};\ x\right)
+(m+1){}_2F_1\left(\begin{array}{c}\alpha,\ \beta\\ \gamma-m\end{array};\ x\right)
\]
and
\begin{eqnarray*}
(\gamma+m)A^\dagger \left({}_2F_1\left(\begin{array}{c}\alpha,\ \beta\\ \gamma+m\end{array};\ x\right)\right)
&=&-(\gamma+m+1-\alpha)(\gamma+m+1-\beta){}_2F_1\left(\begin{array}{c}\alpha,\ \beta\\ \gamma+m+1\end{array};\ x\right)\\
&&\ +m(\gamma+m){}_2F_1\left(\begin{array}{c}\alpha,\ \beta\\ \gamma+m\end{array};\ x\right).
\end{eqnarray*}

Instead of dividing $\mathscr{H}$  by $\partial$ mentioned in Section \ref{SS:pseries}, we consider division by $H$ which gives solutions of the Heun equation $\mathscr{H}f=0$ in $\mathcal{D}/\mathcal{D}H$. This implies that certain solutions 
 of the Heun equation $\mathscr{H}f=0$  can assume the form
\[f(x)=\sum_{m=0}^{\infty}\tilde{C}_m A^{m+1}
\left({}_2F_1\left(\begin{array}{c}\alpha,\ \beta\\ \gamma\end{array};\ x\right)\right)
=\sum_{m=0}^{\infty}{C}_m \cdot {}_2F_1\left(\begin{array}{c}\alpha,\ \beta\\ \gamma-m\end{array};\ x\right).
\]
We note that ${}_2F_1\big(\alpha,\ \beta;  \gamma;\ x\big)$ can be replaced by any other local solutions of $H=0$. 
\bigskip

\begin{remark} The finite sum of the hypergeometric functions
\[\sum_{k=0}^{M}{C}_k \cdot {}_2F_1\left(\begin{array}{c}\alpha,\ \beta\\ \gamma-k\end{array};\ x\right)\]
was also considered in \cite[Shiga, Tsutsui, Worfart]{STW} and
viewed as the solution of the Fuchsian equation with the Riemann $P$-scheme \[
		P_{\mathbb{CP}^1}
		\begin{Bmatrix}
\ 0\ &\ 1\ &\ a_1\ &\ \cdots &\ a_M & \infty\ &\\
0&0&0&\cdots&0&\alpha&x \\
1-\gamma&1-\delta&2&\cdots&2&\beta&
	\end{Bmatrix},
	\]
	where $a_1,\cdots,a_M$ are all apparent singularities.
\end{remark}

Similarly, we define the ``annihilation'' operator $\hat{A}$ and ``creation'' operator $\hat{A}^\dagger$ which are defined by
	\begin{equation}\label{E:A2-Heun}
	\hat{A}=x\partial+\gamma+\epsilon \quad \mbox{ and }\quad \hat{A}^\dagger=(x-1)\partial+\delta-1.
	\end{equation}

In this case, we divide ${\mathscr{H}}$ by $\hat{H}$
in $\mathbb{C}[x,\, \hat{A}^\dagger]$ which gives solutions
of Heun equation $\mathscr{H}f=0$ as the following series
\[f(x)=\sum_{m=0}^{\infty}\tilde{C}_m (\hat{A}^\dagger)^{m}
\left({}_2F_1\left(\begin{array}{c}\alpha,\ \beta\\ \gamma+\epsilon\end{array};\ x\right)\right)
=\sum_{m=0}^{\infty}{C}_m \cdot {}_2F_1\left(\begin{array}{c}\alpha,\ \beta\\ \gamma+\epsilon+m\end{array};\ x\right).\]
\\

We now give the details of the above division algorithms. 
\bigskip

\begin{lemma}\label{lem1} Let us assume the notations \eqref{E:heun}, \eqref{E:H1-Heun}, \eqref{E:H2-Heun}, \eqref{E:A1-Heun}, \eqref{E:A2-Heun} introduced earlier. Then 
\begin{enumerate}
\item[(i)] \begin{eqnarray*}
                 \mathscr{H}&=&(x-a)\mathcal{H}+\epsilon (a-1)X\partial+\alpha\beta a=(X-a)\mathcal{H}+\epsilon (a-1)A+c_1\\
                              &=&(X-a)\hat{\mathcal{H}}+\epsilon a(X-1)\partial+\alpha\beta a=(X-a)\hat{\mathcal{H}}+\epsilon a\hat{A}^\dagger+\hat{c}_1,
           \end{eqnarray*}
           where $c_1=\alpha\beta a-\gamma\epsilon(a-1)$ and $\hat{c}_1=\alpha\beta a-(\delta-1)\epsilon a$.
\item[(ii)] 
			\[
				\begin{split}
					\mathcal{H} &=AA^\dagger+c_2\\
					\hat{\mathcal{H}}&=\hat{A}\hat{A}^\dagger+\hat{c}_2
				\end{split}
			\] where $c_2 =\alpha\beta+\gamma(1-\delta-\epsilon)$ and $ \hat{c}_2 =\alpha\beta+(\gamma+\epsilon)(1-\delta)$;                  
\item[(iii)] 
		\[
			\begin{split}
				\partial &=[A,A^\dagger]=[\partial,A^\dagger]=[\partial,A]\\
				&=[\hat{A},\hat{A}^\dagger]=[\partial,\hat{A}^\dagger]=[\partial,\hat{A}]\\
				&                      ={A}-{A}^\dagger+{c}_3\\
                           &=\hat{A}-\hat{A}^\dagger+\hat{c}_3,
                           \end{split}
                   \]where
      $c_3=\delta+\epsilon-\gamma-1$ 
			and $\hat{c}_3=\delta-\epsilon-\gamma-1$.
\item[(iv)] \[
			\begin{split}
				[\mathcal{H},A] &=-A\partial=\mathcal{H}-A^2-c_3A-c_2\\
				[\hat{\mathcal{H}},\hat{A}^\dagger]
                           & =\partial \hat{A}^\dagger=\hat{\mathcal{H}}-(\hat{A}^\dagger)^2+\hat{c}_3 \hat{A}^\dagger-\hat{c}_2;
			\end{split}
		\]
\item[(v)] \[
			\begin{split}
				&[A,\, X]=X=[\hat{A},\, X]\\
				&[A^\dagger,\, X-1]=X-1=[\hat{A}^\dagger,\, X-1].
			\end{split}
		\]
\item[(vi)] \[
			X(A^2+c_3A+c_2)=X\mathcal{H}+(A-1)(A-\gamma).
		\]
\end{enumerate}
\end{lemma}
\begin{proof}
Routine computation.
\end{proof}
\bigskip

\begin{remark} We note that part (ii) of the above lemma  consists of what is called  ``factorisation formalae" and they can be found in \cite[\S3.4]{Derezinski2014}. On the other hand, part (iv) of the above lemma is a consequence of parts (ii) and (iii).  It can be found in \cite[\S3.4]{Derezinski2014} and  \cite{Der_Maj}. These are called  \textit{transmutation formulae} which play an important role in the theory of special functions. 
\end{remark}
\bigskip

Now we have the following division algorithm
\begin{proposition}\label{prop1}
For each $S\in\mathbb{C}[A]$ with $\deg s=n$, there exist $P_S,Q_S,R_S\in\mathbb{C}[A]$ with
$\deg P_S=\deg Q_S=n$ and $\deg R_S\leq n+1$
such that
	\[
		\mathscr{H} S= (P_SX+Q_S)\mathcal{H}+R_S.
	\]
\end{proposition}
\begin{proof}[Sketch of Proof] By linearity, we suffice to show the statement is true when $S=A^n$ for any $n=0,1,2,\cdots$.
We prove it by induction. First, we can see from part (i) of Lemma \ref{lem1} that
\[
	\mathscr{H}1=(P_1 X+ Q_1)\mathcal{H}+R_1,
\]where $P_1=1$, $Q_1=-a$ and $R_1=\epsilon (a-1)A+c_1$.
Now we assume that
\[\mathscr{H} A^k = (P_{A^k} X+ Q_{A^k})\mathcal{H}+R_{A^k}\]
for some $P_{A^k},\, Q_{A^k},\, R_{A^k}\in\mathbb{C}[A]$ with
$\deg P_{A^k}=\deg Q_{A^k}=k$ and $\deg R_{A^k}\leq k+1$. Then it follows from
Lemma~\ref{lem1} that
		\[
			\begin{split}
				\mathscr{H} A^{k+1}& =[(A-1)P_{A^k}X+(A+1)Q_{A^k}]\mathcal{H}\\
					&\quad -(A-\gamma)(A-1)P_{A^k}
                        -(A^2+c_3A+c_2)Q_{A^k}+AR_{A^k}.
			\end{split}
		\]In particular, the coefficient of $P_{A^k}$ above follows from part (vi) of Lemma \ref{lem1}. 
So the proof follows by taking $P_{A^{k+1}}=(A-1)P_{A^k}$, $Q_{A^{k+1}}=(A+1)Q_{A^k}$ and
		\begin{equation}\label{E:remainder_1}
			R_{A^{k+1}}=-(A-\gamma)(A-1)P_{A^k}-(A^2+c_3A+c_2)Q_{A^k}+AR_{A^k}
		\end{equation}which has degree $k+2$ at most. Thus this completes the induction.
		
\end{proof}
\bigskip

\begin{lemma}\label{R:takemura} From the recursion \eqref{E:remainder_1} and those similar to it, we have, for each $n\in\mathbb{N}\cup\{0\}$,
	\begin{enumerate}
		\item $P_{A^n}=(A-1)^n$,
		\item $Q_{A^n}=-a(A+1)^n$,
		\item the leading coefficient of $R_{A^n}$ is $(\epsilon+n)(a-1)$.
	\end{enumerate}
\end{lemma}
\bigskip

\begin{proof} Direct verification. \end{proof}
\bigskip

\begin{definition}
Let the remainder map $\Phi:\mathbb{C}[A]\rightarrow\mathbb{C}[A]$ be defined by $\Phi(S)=R_S$ for $s\in\mathbb{C}[A]$. We will also endow $\mathbb{C}[A]$ 
with the $A$-adic valuation $\nu$.
\end{definition}

\bigskip

Similarly, we have the division algorthm in $\mathbb{C}[x,\, \hat{A}^\dagger]$
with the dividend~$\mathscr{H}$ and the divisor $\hat{\mathcal{H}}$:
\begin{proposition}
For each $S\in\mathbb{C}[\hat{A}^\dagger]$ with $\deg S=n$, there exist $\hat{P}_S,\, \hat{Q}_S,\, \hat{R}_S\in\mathbb{C}[\hat{A}^\dagger]$ with
$\deg \hat{P}_S=\deg \hat{Q}_S=n$ and $\deg \hat{R}_S\leq n+1$
such that
\[\mathscr{H} S = [\hat{P}_S(X-1)+ \hat{Q}_S]\hat{\mathcal{H}}+\hat{R}_S.\]
\end{proposition}
\begin{proof} By induction (with the similar argument in the proof of Proposition \ref{prop1}).
\end{proof}

\begin{remark} In fact, we have $\hat{P}_{(\hat{A}^\dagger)^n}=(\hat{A}^\dagger-1)^n$ and $\hat{Q}_{(\hat{A}^\dagger)^n}=(1-a)(\hat{A}^\dagger+1)^n$. Moreover,
The leading coefficient of $\hat{R}_{(\hat{A}^\dagger)^n}$ is $(\epsilon+n)a$.
\end{remark}

\begin{definition}
Let $\hat\Phi:\mathbb{C}[\hat{A}^\dagger]\rightarrow\mathbb{C}[\hat{A}^\dagger]$ be defined by $\hat\Phi(S)=R_S$ for $S\in\mathbb{C}[\hat{A}^\dagger]$.
\end{definition}
\bigskip

\section{Finite Dimensional Invariant Subspace 
         and Corresponding Eigensolutions} \label{eigen}

\subsection{Polynomial Solutions}
Let us suppose that we choose either $\alpha=-n$ or $\beta=-n$ ($n=0,1,2,\cdots$) in the Heun equation \eqref{E:heun}. Then for any $S\in\mathbb{C}[X]$ with $\deg S\leq n$,
$\deg{\Psi}_S\leq n$ where the $\Psi_s$ denotes the remainder map \eqref{E:heun_remainder}.
So the space of polynomials of degree at most $n$, denoted by $\mathcal{P}_n[X]$, is invariant  under ${\Psi}$.

So in this case, finding a polynomial solution is equivalent to solving the eigenvalue problem:
\[{\Psi}S^*=q^*S^*\] for some $q^*\in\mathbb{C}$ and non-trivial $S^*\in\mathcal{P}_n[X]$.

\subsection{Factorization of generalized hypergeometric operators}\label{SS:factor_Heun}

We apply the idea of finding an invariant subspace for remainder map encountered above to solve the factorization problem below.

\begin{definition}\label{D:gen_hyp_oper} Let the {\it generalized hypergeometric operator} $L_{a_1,a_2,\cdots,a_p;\ b_1,b_2\cdots,b_q}$
                   to be	
\[L_{a_1,a_2,\cdots,a_p;\ b_1,b_2\cdots,b_q}:=
\partial(X\partial+b_1-1)\cdots(X\partial+b_q-1)-(X\partial+a_1)\cdots(X\partial+a_p).\]
\end{definition}
\bigskip

We note that 

\[{}_{p}F_{q}\left(\begin{array}{c}a_1,\ a_2,\ \cdots ,\ a_p\\ b_1,\ \cdots ,\ b_q\end{array};\ x\right)=
  \sum_{n=0}^\infty\frac{(a_1)_n(a_2)_n\cdots(a_p)_n}{(b_1)_n(b_2)_n\cdots(b_q)_nn!}x^n\]
is the only solution to $L_{a_1,a_2,\cdots,a_p;\ b_1,b_2\cdots,b_q}$ which is analytic and non-vanishing at the origin.
\bigskip

We mention that \cite{LVW} study the lowest order of the differential operator satisfied by the generalised hypergeometric function
\[
	{}_{p+r}F_{q+r}\left(\begin{array}{c}a_1,\ a_2,\ \cdots ,\ a_p,\, e_1+1,\cdots, e_r+1\\ b_1,\ \cdots ,\ b_q,\, e_1,\cdots, e_r\end{array};\ x\right).
\]

Thus, we have the following

\begin{theorem}\label{Take} Let $\mathscr{H}$ be a Heun operator with $\epsilon=-n\in-\mathbb{N}$. Then there exist $q^\ast$, $e_1,\, \cdots. e_n\in \mathbb{C}$ such that
\[L_{\alpha,\beta,e_1+1,\cdots,e_n+1;\ \gamma,e_1,\cdots,e_n}=Q(\mathscr{H}-q^*)\]
for some $Q\in C(X)[\partial]$ with order  $n$.
\end{theorem}
\bigskip

\begin{proof} 
When $\epsilon=-n$ ($n=0,1,2,\cdots$), then the Remark \ref{R:takemura} guarantees that for any $S\in\mathbb{C}[A]$ with $\deg S=n$, $\deg{\Phi}_S\leq n$.
So the space of polynomials of degree at most $n$, denoted by $\mathcal{P}_n[A]$, is  invariant under ${\Phi}$.

Thus in this case, finding a solution which is a finite sum of hypergeometric functions (that is the case that the point $a$ is an apparent singularity)
is equivalent to solving the eigenvalue problem:
\[{\Phi}S^*=q^*S^*\] for some $q^*\in\mathbb{C}$ and non-trivial $S^*\in\mathcal{P}_n[A]$.
In other words, we have
	\begin{equation}
		\mathscr{H}\cdot S^*=Q_1\mathcal{H}+q^\ast\, S^\ast,
	\end{equation}
	for some $Q_1\in \mathbb{C}[X,\, \partial]$ is a certain differential operator and $H$ is a hypergeometric operator. That is, we have
	\begin{equation}\label{E:variant}
		(\mathscr{H}-q^*)\cdot S^*\ _2F_1\left(\begin{array}{c}\alpha,\ \beta\\ \gamma\end{array};\ x\right)=0.
	\end{equation}

Now we can factorize $S^*\in\mathcal{P}_n[A]$ into linear factors
\[ \left(\frac{X}{e_1}\partial+1\right)\left(\frac{X}{e_2}\partial+1\right)\cdots\left(\frac{X}{e_n}\partial+1\right).\]

However, it is well-known and also straightforward to verify that
\[S^*\ _2F_1\left(\begin{array}{c}\alpha,\ \beta\\ \gamma\end{array};\ x\right)=
\ {}_{n+2}F_{n+1}\left(\begin{array}{c}\alpha,\ \beta,\ e_1+1,\ \cdots ,\ e_n+1\\ \gamma,\ e_1,\ \cdots ,\ e_n\end{array};\ x\right)\]
holds, implying that $S^\ast{}_2F_1(\alpha,\, \beta;\, \gamma;\, x)$ is also a solution of $L_{\alpha,\beta,e_1+1,\cdots,e_n+1;\ \gamma,e_1,\cdots,e_n}=0$. To prove that every solution of $\mathscr{H}-q^\ast=0$ is also a solution of $L_{\alpha,\beta,e_1+1,\cdots,e_n+1;\ \gamma,e_1,\cdots,e_n}=0$, we only need to observe that $S^\ast$ of any linear combination of two linearly independent solutions of $H$, also satisfies \eqref{E:variant} and $L_{\alpha,\beta,e_1+1,\cdots,e_n+1;\ \gamma,e_1,\cdots,e_n}=0$ after a routine calculation.
\end{proof}
\bigskip

The above theorem is conjectured  by Takemura in \cite{Take4} for the case of Heun equation having only one apparent singularity.

Thus, by the above theorem, we have
\begin{corollary}[(Heun Identity)]\label{C:identity} Let $\mathscr{H}$ be a Heun operator with $\epsilon=-n\in-\mathbb{N}$. Then there exist a $q^\ast$, $e_1,\, \cdots. e_n\in \mathbb{C}$ such that the local Heun function
	\begin{equation}\label{E:guage}
		\begin{split}
			Hl(a,q^*;\alpha,\beta,\gamma,\delta,\epsilon=-n;\ x) &=
				{}_{n+2}F_{n+1}\left(\begin{array}{c}\alpha,\ \beta,\ e_1+1,\ \cdots ,\ e_n+1\\ \gamma,\ e_1,\ \cdots ,\ e_n\end{array};\ x\right)\\
				&=\left(\frac{X}{e_1}\partial+1\right)\left(\frac{X}{e_2}\partial+1\right)\cdots\left(\frac{X}{e_n}\partial+1\right){}_2F_1\left(\begin{array}{c}\alpha,\ \beta\\ \gamma\end{array};\ x\right).
		\end{split}
	\end{equation}	
\end{corollary}
\bigskip

The above formula is the generalization of the identities given by Maier \cite{Maier}:
\begin{example}[(See Maier \cite{Maier})]
\begin{eqnarray*} Hl(a,q^*;\alpha,\beta,\gamma,\delta,\epsilon=-1;\ x)
&=&{}_{3}F_{2}\left(\begin{array}{c}\alpha,\ \beta,\ e_1+1\\ \gamma,\ e_1\end{array};\ x\right)\\
&=&\left(\frac{X}{e_1}\partial+1\right){}_2F_1\left(\begin{array}{c}\alpha,\ \beta\\ \gamma\end{array};\ x\right),
\end{eqnarray*}
where $a=\frac{e_1(e_1-\gamma+1)}{(e_1-\alpha)(e_1-\beta)}$ and
$q^*=\frac{\alpha\beta(e_1+1)(e_1-\gamma+1)}{(e_1-\alpha)(e_1-\beta)}.$
\bigskip

Thus, by the symmetry of Heun equations, he  re-derived the following ${}_{3}F_{2}$ identities
\[{}_{3}F_{2}\left(\begin{array}{c}\alpha,\ \beta,\ e_1+1\\ \gamma,\ e_1\end{array};\ x\right)
  =(1-x)^{-\alpha}{}_{3}F_{2}\left(\begin{array}{c}\alpha,\ \gamma-\beta-1,\ e'_1+1\\ \gamma,\ e'_1\end{array};\ \frac{x}{x-1}\right)
	\ \mbox{ (Pfaff-like identity)}\]
and
\[{}_{3}F_{2}\left(\begin{array}{c}\alpha,\ \beta,\ e_1+1\\ \gamma,\ e_1\end{array};\ x\right)
=(1-x)^{\gamma-\alpha-\beta-1}{}_{3}F_{2}\left(\begin{array}{c}\gamma-\alpha-1,\ \gamma-\beta-1,\ e''_1+1\\ \gamma,\ e''_1\end{array};\ x\right)
\ \mbox{ (Euler-like identity)},\]
where
$e'_1=\frac{(\gamma-\beta-1)e_1}{e_1-\beta}$ and $e''_1=\frac{(\gamma-\alpha-1)(\gamma-\beta-1)e_1}{e_1(\gamma-\alpha-\beta-1)+\alpha\beta}$.
\end{example}
\bigskip

We quote the following symmetry results of Heun equations from the work of Maier \cite[pp. 825-826]{Maier1}
\begin{proposition} Let $Hl(a, q; \alpha,\, \beta,\, \gamma,\, \delta,\,; x)$ denote a usual local Heun function. Then
	\begin{enumerate}
		\item 
			\begin{equation}
				Hl(a, q; \alpha,\, \beta,\, \gamma,\, \delta\,; x)
				=(1-x)^{1-\delta} Hl (a,\, q-(\delta-1)\gamma a;\, \beta-\delta+1,\, \alpha-\delta+1,\, \gamma,\, 2-\delta; \ x)
			\end{equation}
		\item 			
			\begin{equation}
				Hl(a, q; \alpha,\, \beta,\, \gamma,\, \delta\,; x)
				=(1-x)^{-\alpha} Hl \big(\frac{a}{a-1},\, \frac{-q+\gamma\alpha a}{a-1};\, \alpha,\, \alpha-\delta+1,\, \gamma,\, \alpha-\beta+1; \ \frac{x}{x-1}\big)
			\end{equation}
			\end{enumerate}
\end{proposition}
\bigskip
		
Applying the Corollary \ref{C:identity} to the above symmetry formulae of the Heun solutions allows us to derive the following identities, in consecutive order, that were derived by Miller and Paris \cite{Miller_Paris_2011}, \cite{Miller_Paris_2013} with different methods.
\bigskip

\begin{corollary}\label{C:identity_3} Let $e_1,\, \cdots, e_r\in \mathbb{C}\backslash\{0\}$, $m_1,\, \cdots, m_r$ are positive integers such that $n=m_1+\cdots+ m_r$. There exist  $f_1,\, \cdots, f_n\in \mathbb{C}\backslash\{0\}$ not necessarily the same in each occurrence below, such that
	\begin{enumerate}
		\item 
			\[
				\begin{split}
				&{}_{r+2}F_{r+1}\left(\begin{array}{c}\alpha,\ \beta,\ e_1+m_1,\ \cdots ,\ e_r+m_r\\ \gamma,\ e_1,\ \cdots ,\ e_r\end{array};\ x\right)\\
			&=(1-x)^{\gamma-\alpha-\beta-n}
			{}_{n+2}F_{n+1}\left(\begin{array}{c}\gamma-\alpha-n,\ \gamma-\beta-n,\ f_1+1,\ \cdots ,\ f_n+1\\ \gamma,\ f_1,\ \cdots ,\ f_n\end{array};\ x\right),
				\end{split}
			\]
		\item \[
				\begin{split}
				&{}_{r+2}F_{r+1}\left(\begin{array}{c}\alpha,\ \beta,\ e_1+m_1,\ \cdots ,\ e_r+m_r\\ \gamma,\ e_1,\ \cdots ,\ e_r\end{array};\ x\right)\\
			&=(1-x)^{-\alpha}
			{}_{n+2}F_{n+1}\left(\begin{array}{c}\alpha,\ \gamma-\beta-n,\ f_1+1,\ \cdots ,\ f_n+1\\ \gamma,\ f_1,\ \cdots ,\ f_n\end{array};\ \frac{x}{x-1}\right),
				\end{split}
			\]
		\end{enumerate}
	Each of the above identities is valid in their respective regions of convergence.
\end{corollary}
\bigskip

\begin{proof} Let us consider the special case when $r=1,\ m_1=2$ above:  
	\[
		\begin{split}
	&{}_{3}F_{2}\left(\begin{array}{c}\alpha,\ \beta,\ e_1+2 \\ \gamma,\ e_1\end{array};\ x\right)\\	
	&=
	{}_{4}F_{3}\left(\begin{array}{c}\alpha,\ \beta,\ e_1+1,\, e_1+2 \\ \gamma,\ e_1,\, e_1+1\end{array};\ x\right)\\
	&=\left(\frac{X}{e_1+1}\partial+1\right)\left(\frac{X}{e_1}\partial+1\right){}_2F_1\left(\begin{array}{c}\alpha,\ \beta\\ \gamma\end{array};\ X\right)\\
	&=Hl(a, q; \alpha,\, \beta,\, \gamma,\, -\gamma+\alpha+\beta+3\,; x)
	\quad \textrm{	for some\ } a\ \textrm{and\ } q\\
	&=(1-x)^{-2+\gamma-\alpha-\beta} Hl(a,\, q+(-2+\gamma-\alpha-\beta)\gamma a;\, -2+\gamma-\alpha,\, -2+\gamma-\beta,\, \gamma,\, -1+\gamma-\alpha-\beta; \ x)\\
	&=(1-x)^{-2+\gamma-\alpha-\beta} {}_{4}F_{3}\left(\begin{array}{c} -2+\gamma-\alpha,\, -2+\gamma-\beta,\ f_1+1,\, f_2+1 \\ \gamma,\ f_1,\, f_2\end{array};\ x\right)\quad \textrm{for some\ } f_1\ \textrm{and}\ f_2.
		\end{split}
	\]
	The general case follows similarly.
\end{proof}

\begin{remark}The origins of the above $f_k$ are from the gauge transformations mentioned in the proof of Theorem \ref{Take}.
\end{remark}

\section{Confluent Heun Equations} \label{S:CHE}

\subsection{Confluent HEs and corresponding Eigensolutions}
Let the confluent Heun operator be
  \[\mathscr{H}c(\alpha,\epsilon,\gamma,\delta):=
     X(X-1)\partial^2+[\gamma(X-1)+\delta X+\epsilon X(X-1)]\partial
     + \alpha\epsilon X,
  \] with the Riemann $P$-scheme \[
		P_{\mathbb{CP}^1}
		\begin{Bmatrix}
\ 0\ &\ 1\ &\ \infty\ &\\
0&0&*&x \\
1-\gamma&1-\delta&*&
	\end{Bmatrix},
		\]
which has an irregular singularity at $\infty$ of rank $1$.				
Let the confluent hypergeometric operator $\mathcal{H}c$ be
\[
	\mathcal{H}c(\alpha,\epsilon,\gamma):=
     X\partial^2+[\gamma+\epsilon X]\partial
     + \alpha\epsilon
\]
with the Riemann $P$-schemes \[
		P_{\mathbb{CP}^1}
		\begin{Bmatrix}
\ 0\ &\ \infty\ &\\
0 & * & -\epsilon x \\
1-\gamma& * &
	\end{Bmatrix}.
	\]
Note that this confluent hypergeometric operator $\mathcal{H}c$ differs from the standard one $K$ in \eqref{E:CH2} by a scaling of the independent variable $\epsilon x$. The purpose of this scaling is to facilitate the long division of $\mathscr{H}c$ by $\mathcal{H}c$.

Recall that $A=X\partial+\gamma$.	It is well-known that
\begin{eqnarray*}\label{e1}
A\left({}_1F_1\left(\begin{array}{c}\alpha\\ \gamma-m\end{array};\ x\right)\right)
=(\gamma-m-1){}_1F_1\left(\begin{array}{c}\alpha\\ \gamma-m-1\end{array};\ x\right)
+(m+1){}_1F_1\left(\begin{array}{c}\alpha\\ \gamma-m\end{array};\ x\right)
\end{eqnarray*}
and
\begin{eqnarray*}
A\left({}_1F_1\left(\begin{array}{c}\alpha+m\\ \gamma\end{array};\ x\right)\right)
=(\alpha+m){}_1F_1\left(\begin{array}{c}\alpha+m+1\\ \gamma\end{array};\ x\right)
+(\gamma-\alpha-m){}_1F_1\left(\begin{array}{c}\alpha+m\\ \gamma\end{array};\ x\right).
\end{eqnarray*}
\bigskip

\begin{lemma}\label{lem2} Let $\mathscr{H}c$ and $\mathcal{H}c$ be the confluent Heun and confluent hypergeometric operators respectively defined above.  Then we have the following
\begin{enumerate}
\item[(i)] 
            $c\mathscr{H}=(X-1)(\mathcal{H}c)+\delta X\partial+\alpha\epsilon=(X-1)(\mathcal{H}c)+\delta A+\alpha\epsilon-\delta\gamma$;
\item[(ii)] $\mathcal{H}c=A(\partial+\epsilon)+(\alpha-\gamma)\epsilon$;
\item[(iii)] $[\partial+\epsilon,\ A]=\partial$;
\item[(iv)] $[\mathcal{H}c,A]=A\partial=\mathcal{H}c-\epsilon A-c_1,$ where $c_1=(\alpha-\gamma)\epsilon;$
\item[(v)] $[A,\, X]=X$.
\end{enumerate}
\end{lemma}
\begin{proof}
Routine computation.
\end{proof}
\bigskip

Now we have the following division algorithm
\begin{proposition}
For each $S\in\mathbb{C}[A]$ with $\deg S=n$, there exist $P_S,Q_S,R_S\in\mathbb{C}[A]$ with
$\deg P_S=\deg Q_S=n$ and $\deg R_S\leq n+1$
such that
	\[(\mathscr{H}c) s = (P_SX+Q_S)(Hc)+R_S.\]
\end{proposition}
\begin{proof} By induction (with the similar argument in the proof of Proposition \ref{prop1}).
\end{proof}	
\bigskip
	
\begin{remark} In fact, we have
$P_{A^n}=(A-1)^n$ and $Q_{A^n}=-(A+1)^n$. Moreover, the leading coefficient of $R_{A^n}$ is $\delta+n$.
\end{remark}
\bigskip

\begin{definition}
Let the remainder map $\Phi:\mathbb{C}[A]\to\mathbb{C}[A]$ be defined so that for each $S\in\mathbb{C}[A]$, $\Phi(S)$ is the unique polynomial in $\mathbb{C}[A]$ such that
\[
	\mathscr{H}_cS=(XP+Q)(\mathcal{H}c)+\Phi(S)\hspace{1cm}\mbox{for some }P,Q\in\mathbb{C}[A].
\]
\end{definition}
\bigskip

\subsection{Factorisation of generalized Confluent Hypergeometric operators}

This section parallels to the \S \ref{SS:factor_Heun}. We have

\begin{theorem} Let  $\mathscr{H}c$ be the confluent Heun operator defined above with $\delta=-n\in-\mathbb{N}$. Then there exists a $q^\ast$, $e_1,\, \cdots, e_n\in\mathbb{C}$  such that
	\[
		L_{\alpha,e_1+1,\cdots,e_n+1;\ \gamma,e_1,\cdots,e_n}=Q(\mathscr{H}c-q^*)
		\]for some $Q\in C(x)[\partial]$ with order $n$, 	where	$L_{\alpha,e_1+1,\cdots,e_n+1;\, \gamma, e_1,\cdots,e_n}$ is the generalized hypergeometric operator defined in \eqref{D:gen_hyp_oper}. 
\end{theorem}
\bigskip

\begin{proof}
Since $\delta=-n$, the space of polynomials of degree at most $n$: $ \mathcal{P}_n[A]$ is invariant under the reminder map $\Phi$. Therefore 
we can consider the eigenvalue problem
\[{\Phi}S^*=q^*S^*\] for some $q^*\in\mathbb{C}$ and non-trivial $S^*\in\mathcal{P}_n[A]$.

In other words, we have
\[
	(\mathscr{H}c-q^*)S^*\ _1F_1\left(\begin{array}{c}\alpha \\ \gamma\end{array};\ -\epsilon x\right)=0.
\]

Now we factorize $S^*\in\mathcal{P}_n[A]$ into linear factors
\[ \left(\frac{X}{e_1}\partial+1\right)\left(\frac{X}{e_2}\partial+1\right)\cdots\left(\frac{X}{e_n}\partial+1\right).\]

However, it is straightforward to verify
\[S^*\ _1F_1\left(\begin{array}{c}\alpha \\ \gamma\end{array};\ -\epsilon x\right)=
\ {}_{n+1}F_{n+1}\left(\begin{array}{c}\alpha,\ e_1+1,\ \cdots ,\ e_n+1\\ \gamma,\ e_1,\ \cdots ,\ e_n\end{array};\ -\epsilon x\right)\]

The remaining proof that every solution of $\mathscr{H}_c-q^\ast=0$ is also a solution of $L_{\alpha,e_1+1,\cdots,e_n+1;\ \gamma,e_1,\cdots,e_n}=0$ is similar to that in  Theorem \ref{Take}.
\end{proof}
\bigskip

Thus, we have the following

\begin{corollary}[(Confluent Heun Identity)] The local confluent Heun function
\[
	\mathcal{H}cl(q^*;\alpha,\epsilon,\gamma,\delta=-n;x)=
{}_{n+1}F_{n+1}\left(\begin{array}{c}\alpha,\ e_1+1,\ \cdots ,\ e_n+1\\ \gamma,\ e_1,\ \cdots ,\ e_n\end{array};\ -\epsilon x\right).\]
\end{corollary}
\bigskip

\begin{example} When $\delta=-1$, the matrix representation of the remainder map $\Phi$ with respect to the basis $\{1,A\}$ is
\[\begin{bmatrix}
  \alpha\epsilon-\delta\gamma & \alpha\epsilon-\gamma\epsilon+\gamma\\
	\delta & \alpha\epsilon-\delta\gamma+\epsilon-\gamma-1
  \end{bmatrix}.\]
Thus,	
\begin{eqnarray*} \mathcal{H}cl(q^*;\alpha,\epsilon,\gamma,\delta=-1;x)
&=&{}_{2}F_{2}\left(\begin{array}{c}\alpha,\ e_1+1\\ \gamma,\ e_1\end{array};\ -\epsilon x\right)\\
&=&\left(\frac{X}{e_1}\partial+1\right){}_1F_1\left(\begin{array}{c}\alpha\\ \gamma\end{array};\ -\epsilon x\right),
\end{eqnarray*}
where $e_1=\frac{\delta}{q^*+2\gamma\delta-\alpha\epsilon+\gamma-\epsilon+1}$ and
$q^*$ satisifies \[q^2+(2\delta\gamma+\gamma-2\alpha\epsilon-\epsilon+1)q
+(\delta\gamma-\alpha\epsilon)^2+\alpha\epsilon(\epsilon-1)+\delta\gamma^2.\]
\end{example}
\bigskip

We now quote the following formula from \cite[Prop. 2.1]{DDMRR1978}
\begin{proposition} Let $Hcl(q;\, \alpha,\, \gamma,\, \delta,\, \epsilon;\, x)$ be a local confluent Heun solution. Then
	\[
		\mathcal{H}cl(q;\, \alpha,\, \gamma,\, \delta,\, \epsilon;\, x)
		=e^{\epsilon x}\mathcal{H}cl(q+\epsilon\gamma;\, \delta+\gamma-\alpha,\, \gamma,\, \delta,\, -\epsilon;\, x).
	\]
\end{proposition}
\bigskip

The above formula obviously reduces to the well-known Kummer identity for ${}_1F_1$ after choosing $\delta=0$ and $\epsilon=-1$:
	\[
		{}_{1}F_{1}\left(\begin{array}{c}\alpha\\ 
		\gamma\end{array};\ x\right)
		=e^x {}_{1}F_{1}\left(\begin{array}{c}\gamma-\alpha\\ 
		\gamma\end{array};\ -x\right).
	\]
\bigskip

The above proposition   leads to  the same statement obtained by Miller and Paris \cite{Miller_Paris_2013} below.

\begin{corollary}\label{C:identity_2} Let $e_1,\, \cdots, e_r\in \mathbb{C}\backslash\{0\}$, $m_1,\, \cdots, m_r$ are positive integers with $n=m_1+\cdots+ m_r$. Then there exist  $f_1,\, \cdots, f_n\in \mathbb{C}\backslash\{0\}$ such that
			\[
				\begin{split}
				&{}_{r+1}F_{r+1}\left(\begin{array}{c}\alpha,\ e_1+m_1,\ \cdots ,\ e_r+m_r\\ \gamma,\ e_1,\ \cdots ,\ e_r\end{array};\ -\epsilon x\right)\\
			&=e^{-\epsilon x}
			{}_{n+1}F_{n+1}\left(\begin{array}{c}\gamma-\alpha-n,\ f_1+1,\ \cdots ,\ f_n+1\\ \gamma,\ f_1,\ \cdots ,\ f_n\end{array};\ \epsilon x\right).
				\end{split}
			\]
\end{corollary}
\bigskip

\begin{proof} The proof follows the same idea as that of Corollary \ref{C:identity_3}.
\end{proof}
\bigskip

\section{Concluding Remarks}\label{S:CR}

We have proposed a framework of formal series converging in appropriate $\mathcal{I}-$adic topology 
to solve Weyl-algebraic linear $\partial-$equations that contains the classical series method to solve linear differential equations as a special case when we take $\partial=\frac{d}{dx}$. The setup can also be applied to solve linear difference equations when $\partial$ is interpreted as difference operators such as the forward difference operator $\partial =\Delta$. We have included how this theory is done to the difference Bessel equations below. A full-fledged theory applied to linear difference equations will appear in separate papers. The advantage of the proposed framework includes concepts that generalise the notions of ordinary points, regular singularities, indicial or characteristic equations, etc, which are classically defined for linear differential equations in the complex variable $x$, apply to linear $\partial-$equations with different valid identifications of $\partial$. 

 The main idea of obtaining an $\mathcal{I}-$adic expansion solution to a given Weyl-algebraic linear $\partial-$equation 
$L$ is a comparison principle that  \textit{compares} solutions $L$ and solutions of another \textit{properly selected} Weyl-algebraic $\partial-$adic equations such that the solutions to the latter equation is better known. Indeed, the classical power series and Frobenius series solution to differential equations with ordinary and regular singular point respectively, can be viewed as a comparison between the power series to solution to $X\partial-n=0$ and Frobenius series and solutions to first order equation $X\partial -\nu=0$ respectively.  The comparison principle between two differential operators is achieved by a generalised division algorithm for these differential operators and a root finding algorithm based on Hensel's lemma, a generalised Newton's method established by Kuhlmann from an ultrametric space viewpoint. 
The infinite (Frobenius) series solutions becomes a consequence of proving the existence of zeros of the remainder maps by the Hensel lemma. In order to have a correct fit of $D-$modules structures  of two differential operators, it is crucial to have a good  choice of maximal ideal $I$, for example,  generated by combination of symbols $X$ and $X\partial$, in the case of Parabolic cylindrical functions expansion  solutions for the double confluent Heun equation (also known as \textit{rotating harmonic oscillators} from the work of Schr\"odinger) 
 Thus although  there can appear many very different kind of special function expansions solutions for higher differential equations from mathematical physics on the surface, they are really the same $X-$adic series (Frobenius) series at an  singularity in disguises. 

Our theory not only recovers this classical framework of power series solutions to differential equations, but it also allows us to view various important series expansions of differential equations at the \textit{higher levels}, such as some Heun class differential equations, by orthogonal polynomials/classical special functions 
Weyl-algebraic framework 

Although the main focus of the Part I of this article is on Frobenius series solutions to linear differential equations, the $D-$module approach developed in this article is completely general that it also applies on linear difference equations as shown by the following example.

\subsubsection{Difference Bessel functions}
We illustrate by an example below that the remainder and Newton-Raphson iterations based on Kuhlmann's theory \cite{Kuhl} on immediate maps of valuation theory he developed not only applies to differential equations, but also to linear difference equations (see for example \cite[Chap. 5]{Beals_Wong_2016}). A full-fledged theory will appear in a separate paper. 

We return to the notation defined in \S\ref{S:ultra} where $X\in \mathfrak{X}$ and $[\partial,\, X]=1$. We write the Bessel operator \eqref{E:bessel} in an algebraic form
		\begin{equation}\label{E:Weyl-besse_eqnl}
			L=(X\partial)^2 +(X^2-n^2),
		\end{equation}
		where we refrain from interpreting the $X: \mathbb{C}[x] \to \mathbb{C}[x]$ at this moment. Then the multiplication map in
		\begin{equation}\label{E:D-map-regular-sing}
			\mathcal{D}/\mathcal{D}L\stackrel{\times S}{\longrightarrow}\mathcal{D}/\mathcal{D}(X\partial-n)
	\end{equation}
		that corresponds to \eqref{E:Bessel_eqn} derived can be written in the notation
		\begin{equation}\label{E:weyl-bessel-fn}
			S=:J_n(X)=X^n/2^n \sum_{k=0}^\infty \frac{(-1)^k(X/2)^{2k}}{\Gamma(n+k+1)k!}
		\end{equation}
	as a member in $\overline{\mathbb{C}[X]}$, which is a \textit{solution} to $LJ_\nu(X)=0\mod X\partial-n$.  We will recover the classical Bessel equation \eqref{E:bessel} $Lf(x)=0$ if we interpret  $\partial=\frac{d}{dx}$ and $X$ as defined by $Xf(x)=xf(x)$. If, however, we interpret $X$ and $\partial$  according to
	\begin{equation}\label{E:forward-shift}
		Xf(x)=xf(x-1)
	\end{equation}	and $\partial=\Delta$ where $\partial f(x) =\Delta f(x)=f(x+1)-f(x)$, so that 	
	\[
		X\cdot 1=x,\qquad X^n\cdot 1=x(x-1)\cdots (x-n+1)
	\]for each integer $n\in\mathbb{N}$, then we still have 
	\[
			[ \partial,\, X]=1
	\]and the equation assumes the form
	\begin{equation}\label{E:Delta-bessel}
			x(x-1)\Delta^2 f(x-2)+x\Delta f(x-1)+x(x-1)y(x-2) -n^2f(x)=0
	\end{equation}
	which is precisely the \textit{difference Bessel equation} obtained previously by Bohner and Cuchta \cite{BC_2017}. The ``solution" \eqref{E:Delta-bessel-fn} after being interpreted in this manner becomes 
	\begin{equation}\label{E:Delta-bessel-fn}
		J_n(X)\cdot 1 = J^\Delta_n(x)=\sum_{k=1}^\infty \frac{(-1)^k  x(x-1)\cdots (x-n-2k+1)}{2^{n+2k}(n+k)!\, k!}
	\end{equation}which, apart from a constant multiple, is also the \textit{difference Bessel equation} proposed by Bohner and Cuchta \cite[Theorem 1]{BC_2017} for which the \eqref{E:Delta-bessel-fn} is a solution\footnote{We discuss its convergence in an half-plane of $\mathbb{C}$ in \cite{Chiang_Ching_Lin}.}. However, the fact that the function \eqref{E:Delta-bessel-fn} being a solution to  \eqref{E:Delta-bessel} is intrinsic from the way in which our theory is developed. It is also straightforward to consider difference Bessel functions for $J^\Delta_\nu(x)$ for non-integer $\nu$ and to obtain corresponding formulae by considering the more general map \eqref{E:D-map-regular-sing} with $n$ being replaced by $\nu$. Thus the ${D}-$module approach to Bessel functions allows us to  derive the difference Bessel equations and functions without extra effort. In particular,  the new difference Bessel functions are unified with the classical  Bessel differential equations and functions under the ${D}-$modules framework. Detailed studies concerning ${D}-$modulus approach on \textit{Bessel calculus} will appear in \cite{Chiang_Ching_Lin}. Our theory applies equally well to  difference versions of Hermite functions, generalised Laguerre, hypergeometric equations and functions, etc.  Indeed, this $D-$modulus approach to some of these special functions as well as their $q-$analogues eventually directs us to study the solution spaces of some systems of holonomic PDEs, a path that appears to have been little trodden in recent decades. Applications  to Gevrey series of arithmetic types \cite{Andre_2000}, \cite{Andre_2003}, \cite{Beukers_2006} will be considered in a future project. 
	
\appendix

\section{Takemura's Conjecture}\label{takemura}

Let the Fuchsian operator $H$ be
  \[H=\partial^2+\left(\frac{\gamma}{x}+\frac{\delta}{x-1}-
    \sum_{k=1}^M\frac{m_k}{x-a_k}\right)\partial+\frac{\alpha\beta x^M+c_{M-1}x^{M-1}\cdots+c_0}{x(x-1)(x-a_1)\cdots(x-a_M)}
  \]
	with $\alpha+\beta-\gamma-\delta+\sum_{k=1}^M m_k+1=0$ and the Riemann $P$-scheme \[
		P_{\mathbb{CP}^1}
		\begin{Bmatrix}
\ 0\ &\ 1\ &\ a_1\ &\ \cdots &\ a_M & \infty\ &\\
0&0&0&\cdots&0&\alpha&x \\
1-\gamma&1-\delta&1+m_1&\cdots&1+m_M&\beta&
	\end{Bmatrix}.
	\]

\begin{conjecture}[(Takemura \cite{Take4})]
Let $m_1,\cdots,m_M$ be positive integers. Suppose that $a_1,\cdots,a_M$ are all apparent singularities in $Hy=0$.
Then there exist $e_1,\cdots,e_N\in\mathbb{C}$ ($N=m_1+\cdots+m_M$) such that
\[L_{\alpha,\beta,e_1+1,\cdots,e_N+1;\ \gamma,e_1,\cdots,e_N}=QH\]
for some $Q\in \mathbb{C}(x)[\partial]$ of order $N$.
\end{conjecture}
\begin{remark}
The above conjecture has been proved by Takemura (involving the calculation by Maple) only for the case: 
\begin{enumerate}
\item $M=1$ together with $m_1=1,2,3,4$ or $5$;
\item $M=2$ together with  $m_1+m_2=1,2,3$ or $4$;
\item $M=3$ together with $m_1=m_2=m_3=1$.
\end{enumerate}
\end{remark}
\begin{remark}
If $H$ is irreducible and $M\in\mathbb{N}$ together with $m_1=m_2=\cdots=m_M=1$, Shiga, Tsutsui and Worfart \cite{STW}
show that the solution of $Hy=0$ is the finite sum of hypergeometric functions
\[\sum_{m=0}^{M}C_m\cdot {}_2F_1\left(\begin{array}{c}\alpha,\ \beta\\ \gamma-m\end{array};\ x\right),\]
and hence the conjecture will also follows for this case.

\end{remark}


\begin{acknowledgements} The authors acknowledge useful conversations with their colleague Henry Cheng.
\end{acknowledgements}

\bibliography{references}
\bibliographystyle{abbrv}


\end{document}